\documentclass[11pt]{article}

\usepackage{graphicx,times,latexsym,amssymb,amsbsy,amsmath,epsfig}

\input xy

\xyoption{all}

\newtheorem{thm}{Theorem}[section]
\newtheorem{prop}[thm]{Proposition}
\newtheorem{lem}[thm]{Lemma}
\newtheorem{cor}[thm]{Corollary}
\newtheorem{dfn}[thm]{Definition}

\renewcommand{\rm}[1]{\mathrm{#1}}
\renewcommand{\cal}[1]{\mathcal{#1}}

\newcommand{\bbR}{\mathbb{R}}

\newcommand{\rmH}{\mathrm{H}}

\newcommand{\m}{\mathrm{m}}

\newcommand{\id}{\mathrm{id}}

\newcommand{\B}{\mathcal{B}}
\newcommand{\C}{\mathcal{C}}

\newcommand{\F}{\mathcal{F}}

\renewcommand{\P}{\mathcal{P}}
\newcommand{\Q}{\mathcal{Q}}

\newcommand{\U}{\mathcal{U}}

\newcommand{\Z}{\mathcal{Z}}

\newcommand{\G}{\Gamma}

\renewcommand{\a}{\alpha}
\renewcommand{\b}{\beta}
\newcommand{\eps}{\varepsilon}
\newcommand{\g}{\gamma}
\renewcommand{\k}{\kappa}
\renewcommand{\l}{\lambda}

\newcommand{\s}{\sigma}
\renewcommand{\phi}{\varphi}

\newcommand{\img}{\mathrm{img}\,}

\renewcommand{\hat}[1]{\widehat{#1}}
\newcommand{\ol}[1]{\overline{#1}}

\newcommand{\into}{\hookrightarrow}
\newcommand{\onto}{\twoheadrightarrow}
\newcommand{\qed}{\nolinebreak\hspace{\stretch{1}}$\Box$}
\newcommand{\fin}{\nolinebreak\hspace{\stretch{1}}$\lhd$}

\renewcommand{\to}{\longrightarrow}

\newcommand{\Seg}{\mathrm{Seg}}
\newcommand{\cts}{\mathrm{cts}}

\newcommand{\al}{\mathrm{al}}
\renewcommand{\sl}{\mathrm{sl}}

\begin{document}

\title{On discontinuities of cocycles in cohomology theories for topological groups}

\author{Tim Austin\footnote{Research supported by a fellowship from the Clay Mathematics Institute}\\ \\ \small{Department of Mathematics, Brown University,}\\ \small{Box 1917, 151 Thayer Street, Providence, RI 02912, U.S.A.}\\ \small{\texttt{timaustin@math.brown.edu}}}

\date{}

\maketitle


\begin{abstract}
This paper studies two cohomology theories for topological groups and modules: Segal's theory based on soft resolutions, and Moore's theory based on a measurable bar resolution.  First it is shown that all classes in Moore's theory have representatives with considerable extra topological structure beyond measurability.  Using similar tools, one can then construct a direct comparison map between Segal's and Moore's theories when both are defined, and show that this map is an isomorphism if the module is discrete.
\end{abstract}

\tableofcontents

\section{Introduction}

Let $G$ be a topological group and $A$ a topological Abelian group on which $G$ acts continuously by automorphisms.  Under a variety of additional assumptions on $G$ and $A$, several proposals have been made for cohomology theories $\rmH^\ast(G,A)$ which parallel the classical cohomology of discrete groups but take the topologies into account.  This paper studies two of these theories, due to Calvin Moore and Graeme Segal.

In~\cite{Moo64(gr-cohomI-II),Moo76(gr-cohomIII),Moo76(gr-cohomIV)}, Moore introduced the theory $\rmH^\ast_\m(G,A)$ based on bar resolutions of measurable cochains.  If $G$ is locally compact and second countable, if one focusses on the category of Polish $G$-modules, and if one requires that `exact sequences' of such modules be algebraically exact, then the resulting theory can be shown to define an effaceable cohomological functor.  It is therefore unique on that category by Buchsbaum's criterion.  It can then be shown to enjoy analogs of all the standard properties of classical group cohomology for these classes of topological groups: for example, if $A$ is also l.c.s.c. then $\rmH_\m^2(G,A)$ classifies topological group extensions of $A$ by $G$.

A more abstract alternative was proposed by Segal in~\cite{Seg70}.  He allows all topological groups $G$ which are groups in the category of k-spaces, and then considers the category of $G$-modules which are Hausdorff k-spaces and are locally contractible.  He also makes the convention that a `short exact sequence' $A\into B\onto C$ must be algebraically exact and must have a local cross-section (that is, $C$ contains an identity neighbourhood on which the quotient map from $B$ has a continuous section).  In this category Segal defines an object to be `soft' it is of the form $\C_\cts(G,A)$ with $A$ a contractible $G$-module, where $\C_\cts$ denotes a space of continuous functions with the compact open topology.  He then shows that any $G$-module in his category admits a rightwards resolution by soft modules, and then that the functor $A\mapsto A^G$ is `derivable' on this category, which means that applying this functor to any choice of soft resolution of $A$ gives a new complex with the same homology.  These homology groups comprise Segal's theory $\rmH^\ast_\Seg(G,A)$, and the standard arguments of homological algebra show that they define a universal cohomological functor on Segal's category of modules for any $G$.

A third theory $\rmH^\ast_{\rm{ss}}$ was introduced by David Wigner in~\cite{Wig73} using semi-simplicial sheaves, and has recently been studied further by Lichtenbaum in~\cite{Lic09} and Flach in~\cite{Fla08}.  It allows any topological group $G$ and $G$-module $A$, and coincides with $\rmH_\Seg^\ast$ for a k-space group and a module in Segal's category, so it is really an extension of Segal's theory.

These different theories have various advantages.  On the one hand, l.c.s.c. groups and Polish modules are the natural setting for most of functional analysis and dynamical systems, and so the universality of $\rmH^\ast_\m$ on that category strongly recommends it for those applications.  However, in other areas, such as class field theory, the sheaf-theoretic definition of $\rmH^\ast_\rm{ss}$ aligns it more closely with cohomologies of other spaces with which it must be compared (see Lichtenbaum's paper for more on this).  Also, the theories $\rmH^\ast_\Seg$ and  $\rmH^\ast_\rm{ss}$ admit spectral sequences that greatly facilitate explicit calculations, and it is not known whether $\rmH^\ast_\m$ can be equipped with any comparable tool.

It is therefore of interest to find cases in which $\rmH^\ast_\m$ and $\rmH^\ast_\Seg$ coincide.  Several cases of agreement have been known for some time, particularly since Wigner's work~\cite{Wig73}. The recent paper~\cite{AusMoo--cohomcty} enlarges the list.  It also contains a much more careful description of how the various theories are defined and the historical context to their study, so the reader is referred there for additional background.  (Those papers also study cases of agreement with another theory, $\rmH^\ast_\rm{cs}$, defined using a classifying space of $G$ and which does not have such obvious universality properties.  That theory is also important for its usefulness in computations, but we will not consider it here.)

For Fr\'echet modules, Theorem A of~\cite{AusMoo--cohomcty} shows that all theories coincide with the theory defined by continuous cochains.  Outside that setting, the strongest comparison results in~\cite{AusMoo--cohomcty} are Theorems E and F.  The heart of these results is the assertion that
\[\rmH^\ast_\m(G,A) \cong \rmH^\ast_\rm{ss}(G,A) \cong \rmH^\ast_\Seg(G,A)\]
whenever $A$ is discrete.  This conclusion is then easily extended to all locally compact and locally contractible $A$ by the Structure Theory for locally compact Abelian groups, an appeal to Theorem A of~\cite{AusMoo--cohomcty} and some diagram-chasing. Note that the second isomorphism here is already clear from the above-mentioned agreement of $\rmH^\ast_\rm{ss}$ and $\rmH^\ast_\Seg$ on Segal's category of modules.

The proof of Theorem F in~\cite{AusMoo--cohomcty} requires several steps.  It relies crucially on breaking up a general group $G$ into its identity component $G_0$ and the quotient $G/G_0$, and then on using the structure of $G_0$ as a compact-by-Lie group promised by the Gleason-Montgomery-Zippin Theorem.  These various special cases are sown together using the Lyndon-Hochschild-Serre spectral sequences for $\rmH^\ast_\m$ and $\rmH^\ast_\rm{ss}$.

In using a separation of cases based on such heavy machinery, an intuitive understanding of why $\rmH^\ast_\m$ and $\rmH^\ast_\rm{ss}$ should agree (in spite of their very different definitions) becomes obscured.  The present paper provides and alternative, more direct proof in case $A$ is discrete.  In that setting we may work with the simpler theory $\rmH^\ast_\Seg$ in place of $\rmH^\ast_\rm{ss}$.

\vspace{7pt}

\noindent\textbf{Theorem A}\quad \emph{If $G$ is an l.c.s.c. group and $A$ is a discrete $G$-module then one has an isomorphism of cohomology theories
\[\rmH^\ast_\m(G,A) \cong \rmH^\ast_\Seg(G,A).\] }

\vspace{7pt}

Owing to the relations that were already known among $\rmH^\ast_\Seg$, $\rmH^\ast_\rm{ss}$ and $\rmH^\ast_\rm{cs}$ prior to the appearance of~\cite{AusMoo--cohomcty}, this essentially recovers the new comparison results of that paper.  Unlike in~\cite{AusMoo--cohomcty}, where $\rmH^\ast_\Seg$ was discussed mostly as a digression, here it will be the fulcrum of this comparison.

To prove Theorem A, we will first introduce two new cohomology theories, denoted $\rmH^\ast_{\rm{sl}}$ and $\rmH^\ast_{\rm{al}}$, which are defined using resolutions consisting of cocycles that have some special topological structure: they are `semi-layered' or `almost layered' functions, respectively.  These notions will be defined in Sections~\ref{sec:sl} and~\ref{sec:al}.  We will then show that one always has $\rmH^\ast_{\rm{sl}}\cong \rmH^\ast_\Seg$ and $\rmH^\ast_{\rm{al}}\cong \rmH^\ast_\m$, and finally observe that in case the target module is discrete it is obvious that $\rmH^\ast_\rm{sl}$ and $\rmH^\ast_\rm{al}$ coincide.  The proofs of these isomorphisms of theories will be fairly simple outings for Buchsbaum's criterion, once the necessary topological preliminaries have been completed.

Importantly, the new theories $\rmH^\ast_\sl$ and $\rmH^\ast_\al$ must be introduced on the same categories of modules as $\rmH^\ast_\Seg$ and $\rmH^\ast_\m$, respectively -- it will not suffice to define them only for discrete modules, say.  This is because if we begin with a discrete module, the induction by dimension-shifting that underlies Buchsbaum's criterion usually converts it into a non-discrete one.  Thus, the formulation of the special classes of cocycles (`semi-layered' and `almost layered') that give rise to $\rmH^\ast_\sl$ and $\rmH^\ast_\al$ can be viewed as the formulation of a successful inductive hypothesis.  It is the main innovation of the present paper.

In the case of Segal's cohomology, his original paper~\cite{Seg70} implicitly offers a concrete resolution for its computation, but of a rather complicated form.  It is a sequence of modules constructed by alternately forming function-spaces and quotients; the ingredients are similar to a bar resolution, but are arranged more intricately.  The resolution that underlies $\rmH^\ast_\sl$ is not `soft' in Segal's sense, but it does give a representation of the same cohomology theory that is closer to the classical bar resolution.  (In this connexion, a recent work of Fuchssteiner, Wagemann and Wockel has provided another such representation.  Our cocycles are quite different from theirs, and can be used in different ways, but we will offer some comparison of these representations later in the paper.)

The methods used to prove Theorem A can also give more elementary results about the usual measurable homogeneous bar resolution, to the effect that all classes have representatives with some additional structure.  The following have some independent interest.

\vspace{7pt}

\noindent\textbf{Theorem B}\quad \emph{If $G$ is an l.c.s.c. group and $A$ is a Polish $G$-module, then any class in $\rmH^p_\m(G,A)$ has a representative cocycle in the homogeneous bar resolution that is continuous on a dense G$_\delta$-set of full measure, including at the origin of $G^{p+1}$.}

\vspace{7pt}

\noindent\textbf{Theorem C}\quad \emph{If $G$ is an l.c.s.c. group and $A$ is a discrete $G$-module, then any class in $\rmH^p_\m(G,A)$ has a representative cocycle in the homogeneous bar resolution that is locally finite-valued and is locally constant on a dense open set of full measure.}

\emph{Moreover, if $G$ is a closed algebraic subgroup of $\rm{GL}_n(\bbR)$ for some $n$ and $A$ is a discrete $G$-module, then a representative $\s$ may be found which is measurable with respect to a partition of $G^p$ into semi-algebraic sets (with reference to the structure of $G^p$ as a real algebraic variety in the real affine space $M_{n\times n}(\bbR)^p$ of $p$-tuples of matrices), and is locally constant at the origin of $G^{p+1}$.}

\vspace{7pt}

\noindent\textbf{Remark}\quad By the usual formula relating cocycles in the homogeneous and inhomogeneous bar resolutions it follows easily that Theorems B and C hold in the latter resolution as well. \fin

\vspace{7pt}

Like Theorem A, the core of Theorems B and C is the formulation of a class of maps from l.c.s.c. groups to Polish modules which all have the properties asserted in those theorems; which include all crossed homomorphisms; and which can be lifted through continuous epimorphisms of target modules and so can be carried to higher degrees by dimension-shifting.  The properties of the cocycles promised by Theorems B and C do not themselves define such a class, so some refinement is necessary, but it turns out that a suitable formulation is rather simpler here than in the case of Theorem A.  We shall therefore prove Theorems B and C first, in Section~\ref{sec:B-and-C}, before formulating further new classes of functions and then using them to complete the proof of Theorem A in Sections~\ref{sec:diss-and-lay} through~\ref{sec:A}.

As the present paper neared completion, my attention was drawn by Christoph Wockel to the preprints~\cite{Fuc11a, Fuc11b,FucWoc11,WagWoc11}.  Those papers explore a variety of cohomology theories for topological groups and modules, including the theory that results from a bar resolution whose cochains are assumed to be continuous on some neighbourhood around the identity, but not globally.  A key theorem of~\cite{WagWoc11} (building on technical results of those other works) asserts that this locally-continuous-cochains theory agrees with $\rmH^\ast_\Seg$ when both are defined.  Knowing this, one can easily construct a comparison map $\rmH^\ast_\Seg(G,A)\to \rmH^\ast_\m(G,A)$ when both theories are defined and then use our Theorem B to show that it is surjective when $A$ is discrete.  However, it still seems tricky to prove injectivity, and hence isomorphism, without something like our more delicate proof of Theorem A below.  We sketch this relation at the end of Section~\ref{sec:B-and-C}.

\subsubsection*{Acknowledgements}

This paper emerged at a tangent from the joint work~\cite{AusMoo--cohomcty}; it benefited greatly from several communications with Calvin Moore.  Christian Wockel made helpful suggestions and brought his paper~\cite{WagWoc11} to my attention.

\section{Preliminaries}\label{subs:prelim}

\subsection*{Basic conventions}

Let $I := (0,1]$ and let $\l$ be Lebesgue measure on $I$.

All topological spaces in this paper will be paracompact; by a theorem of Stone this includes all metrizable spaces (see, for instance, M.E. Rudin~\cite{Rud(M)69}).  The reader will lose little by thinking of all our spaces as Polish.

If $A$ is a Polish Abelian group then we let $LA$ denote the group of $\l$-equivalence classes of measurable functions $I\to A$, and give $LA$ the topology of convergence in measure.  For example, if $A = \bbR$ then $LA = L^0(\bbR)$ with its customary topology.

On the other hand, for any Hausdorff topological Abelian group $A$ we let $EA$ denote the subgroup of left-continuous step functions $I\to A$ with only finitely many discontinuities.  This may be expressed as $\bigcup_{n\geq 1}E^{(n)}A$ with $E^{(n)}A$ the subset of functions having at most $n$ discontinuities.  Unless stated otherwise, we will consider $EA$ as endowed with the direct limit of the topologies on the subsets $E^{(n)}A$, where those topologies are given by the identification of $E^{(n)}A$ with a quotient of $\Delta_n\times A^{n+1}$, where $\Delta_n \subseteq \bbR^{n+1}$ is the $n$-simplex (see~\cite{Seg70}). If $A$ is Polish, this is the topology on $E^{(n)}A$ inherited from $LA$, but the resulting direct limit topology on the whole of $EA$ is usually strictly finer than the topology that $EA$ itself inherits as a subspace of $LA$.

Let $\iota:A\into LA$ or $\iota:A\into EA$ denote the inclusion of $A$ as the constant functions.  The following basic facts are proved by Segal in Proposition A.1 of~\cite{Seg70}.

\begin{prop}\label{prop:Segal-x-sec}
The topological group $EA$ is contractible, and the subgroup $\iota(A)$ has a local cross-section in $EA$. \qed
\end{prop}

\subsection*{Segal cohomology}\label{subs:Seg}

Let $G$ be any topological group in the category of k-spaces, and let $A$ be any topological $G$-module that is likewise a k-space and is locally contractible.  When a choice of $G$ is understood, we will refer to this as \textbf{Segal's category} of modules.  In this category a short exact sequence of continuous module homomorphisms is \textbf{distinguished} if the quotient homomorphism has a local continuous cross-section as a map between topological spaces.

Segal's cohomology for such groups and modules is defined in terms of a fairly abstract class of resolutions.

Such a $G$-module is \textbf{soft} if it takes the form $\C_\cts(G,B)$ for some contractible $G$-module $B$, where this denotes the space of continuous functions $G\to B$ with the compact-open topology and with the diagonal $G$-action. Any $A$ in Segal's category may be embedded into a soft module via the composition of the embeddings
\begin{eqnarray}\label{eq:Seg-cplx}
A\stackrel{\iota}{\to} EA \stackrel{\rm{consts}}{\to} \C_\cts(G,EA) =: E_GA.
\end{eqnarray}
By Proposition~\ref{prop:Segal-x-sec} and the easy fact that $EA$ has a global cross-section in $\C_\cts(G,EA)$ (for instance, by evaluating at $e \in G$), the image of $A$ under this embedding has a local cross-section in $E_GA$.  Forming the quotient module $B_GA := E_GA / A$ therefore gives a distinguished short exact sequence in Segal's category.  Iterating this construction gives a resolution of $A$ by soft modules
\[A \to E_GA \to E_GB_GA \to E_G B_G^2 A \to \ldots\]
(see Proposition 2.1 in~\cite{Seg70}). Now applying the fixed-point functor $A\mapsto A^G$ to this sequence, the resulting homology groups are the \textbf{Segal cohomology groups} $\rmH^\ast_\Seg(G,A)$.

Segal proves in~\cite{Seg70} that this is a universal definition in the sense that any other soft resolution of $A$ gives the same cohomology groups (the fixed-point functor is `derivable', in his terminology).  Importantly, this leads to universality in the sense of Buchsbaum~\cite{Buc60}, in exact analogy with the universality of derived functors in classical homological algebra.  The identity $\rmH^0_\Seg(G,A) = A^G$ and the fact that classes are always effaced under the inclusion $A\into \C_\cts(G,EA)$ are built into Segal's definition, and the existence of long exact sequences follows as an easy exercise (Proposition 2.3 in~\cite{Seg70}).  Therefore, in order to prove  that another candidate theory gives the same cohomology groups as Segal's, one need only check that it has these three properties on Segal's category of modules.

\vspace{7pt}

\noindent\textbf{Remark}\quad Another resolution of $A$ suggested by Segal's theory is 
\[A \to \C_\cts(G,EA)\to \C_\cts(G^2,E^2A)\to \cdots.\]
I do not know whether this is always soft in Segal's sense --- in particular, whether it admits local cross-sections --- and so offers an easier route to calculations in $\rmH^\ast_\Seg$.  This seems unlikely in general, but even if it fails it would be interesting to know more about the homology obtained by applying $(-)^G$ to this resolution. \fin

\subsection*{Measurable cohomology}

We will use the definition of $\rmH^\ast_\m$ based on the measurable homogeneous bar resolution.  As for discrete cohomology, one obtains the same theory from the inhomogeneous bar resolution; this equivalence follows from a routine appeal to Buchsbaum's criterion as in Theorem 2 of~\cite{Moo76(gr-cohomIII)}.

For a l.c.s.c. group $G$, Polish $G$-module $A$ and integer $p\geq 0$ we let $\C(G^p,A)$ denote the group of Haar-a.e. equivalence classes of measurable functions $G^p\to A$, interpreting this as $A$ itself when $p=0$.  This is also a Polish group in the topology of convergence in measure on compact subsets, and if $A$ carries a continuous action of $G$ by automorphisms then we equip each $\C(G^p,A)$ with the associated \textbf{diagonal} action:
\[(g\cdot \phi)(g_1,g_2,\ldots,g_p) = g\cdot\big(\phi(g^{-1}g_1,g^{-1}g_2,\ldots,g^{-1}g_p)\big).\]
We also sometimes write $\C^p(G,A) := \C(G^p,A)$.

With this in mind, one forms the exact resolution of $A$ given by
\[A\stackrel{d}{\to} \C(G,A)\stackrel{d}{\to} \C(G^2,A)\stackrel{d}{\to}\C(G^3,A)\stackrel{d}{\to}\ldots\]
with the usual differentials defined by
\[d\s(g_1,\ldots,g_{p+1}) := \sum_{i=1}^{p+1}(-1)^{p+1 - i}\s(g_1,\ldots,\hat{g_i},\ldots,g_{p+1})\]
for $\s \in \C(G^p,A)$, where the notation $\hat{g_i}$ means that the entry $g_i$ is omitted from the argument of this instance of $\s$.  Note our convention is that the last term always has coefficient $+1$: this avoids some other minus-signs later. Now omitting the initial appearance of $A$ and applying the fixed-point functor $A\mapsto A^G$ gives the complex
\begin{eqnarray}\label{eq:bar-cplx}
\C(G,A)^G\stackrel{d}{\to} \C(G^2,A)^G\stackrel{d}{\to}\C(G^3,A)^G\stackrel{d}{\to}\ldots.
\end{eqnarray}
Letting $\Z^p(G,A) := \ker d|_{\C(G^{p+1},A)^G}$ and $\B^p(G,A) := \img d|_{\C(G^p,A)^G}$, Moore's \textbf{measurable cohomology groups} of the pair $(G,A)$ are the homology groups
\[\rmH^p_\m(G,A) := \frac{\Z^p(G,A)}{\B^p(G,A)}.\]

The basic properties of this theory can be found in~\cite{Moo64(gr-cohomI-II),Moo76(gr-cohomIII),Moo76(gr-cohomIV)}, including the existence of long exact sequences, effaceability, and interpretations of the low-degree groups.  For reference, let us recall that a class in $\rmH^p_\m(G,A)$ may always be effaced using the constant-functions inclusion $A\into \C(G,A)$.  More explicitly, given a cocycle $\s:G^{p+1}\to A$ in the complex~(\ref{eq:bar-cplx}), one has $\s = d\psi$ with $\psi:G^p\to \C(G,A)$ defined by
\begin{eqnarray}\label{eq:basic-dim-shift}
\psi(g_1,\ldots,g_p)(g) := \s(g_1,\ldots,g_p,g)
\end{eqnarray}
(where our choice of signs in the formula for $d$ avoids the need for a minus-sign here).

A theory satisfying all of these properties on the category of Polish $G$-modules is universal by Buchsbaum's criterion, and this fact forms the basis for a comparison with other possible cohomology theories.

\section{Warmup: additional regularity for cocycles}\label{sec:B-and-C}

\subsection*{Proofs of Theorems B and C}

In this section we prove Theorems B and C, which concern only the measurable-cochains theory in the homogeneous bar resolution.  The rest of the paper will go towards proving Theorem A, which requires ideas that are related, but more complicated.  The key point is to define classes of functions that enhance the conclusions of Theorems B and C and which give a hypothesis that can be closed on itself in a dimension-shifting induction.

\begin{dfn}\label{dfn:type-I-and-II}
If $X$ is a locally compact and second countable metrizable space, $\mu$ is a Radon measure of full support on $X$, and $A$ is a Polish Abelian group, then a map $f:X\to A$ is of \textbf{type I} if it is locally finite-valued and there is an open subset $U\subseteq X$ of full $\mu$-measure on which $f$ is locally constant.  It is \textbf{almost type-I} if it is a locally uniform limit of type-I functions.

If, in addition, $X$ is a real algebraic variety with its Euclidean topology and $\mu$ is a smooth measure, then a function $f:X\to A$ is of \textbf{type II} if it takes locally finitely many values and its level sets agree locally with semi-algebraic subsets of $X$.  It is \textbf{almost type-II} if it is a locally uniform limit of type-II functions.

Finally, if $f$ is an almost type-I (resp. almost type-II) function and $x_0 \in X$, then $f$ is \textbf{regular at $x_0$} if it is a limit of type-I (resp. type-II) functions each of which is locally constant around $x_0$ (possibly with different neighbourhoods of constancy).
\end{dfn}

As usual, for locally compact $X$, `locally uniform' convergence refers to convergence in the compact open topology.  In all the cases that follow $X$ will be $G^p$ for some l.c.s.c. group $G$ and $\mu$ will be a left-invariant Haar measure. The basic properties of real algebraic varieties and semi-algebraic sets can be found, for instance, in Bochnak, Coste and Roy~\cite{BocCosRoy98}.  We will not need any sophisticated theory for them here.  It is easy to see that (almost) type-II is stronger than (almost) type-I when both notions make sense. The first simple properties that we need are contained in the following lemmas.

\begin{lem}[Slicing]\label{lem:slice}
If $G$ is an l.c.s.c. group, $m_G$ a left-invariant Haar measure and $f:G^{p+1}\to A$ an almost type-I function, then for almost every $h \in G$ the slice
\[f_h: G^p\to A:(g_1,\ldots,g_p)\mapsto f(g_1,g_2,\ldots,g_p,h)\]
defines an almost type-I function $G^p\to A$.  If $G$ is an algebraic subgroup of $\rm{GL}_n(\bbR)$ then the same holds with `type-II' in place of `type-I'.

If $f$ is equivariant then these properties hold for strictly every $h$, and if $f$ is also regular at the identity then $f_h$ is regular at $(h,h,\ldots,h)$.
\end{lem}

\noindent\textbf{Proof}\quad Let $(\g_n)_n$ be a sequence of type-I (or, where applicable, type-II) functions that converge locally uniformly to $f$. For each $n$, let $U_n$ be a full-measure open set on which $\g_n$ is locally constant.  We need only observe that the intersections
\[(G^p\times \{h\})\cap U_n\]
are all still open, and by Fubini's Theorem they still have full measure for a.e. $h$.  Also, if $G$ is algebraic and $\partial U_n$ is semi-algebraic, then so are these intersections.  Hence for a.e. $h$ the restrictions
\[(g_1,\ldots,g_p)\mapsto \g_n(g_1,g_2,\ldots,g_p,h)\]
are still of type I (or, where applicable, type II), and $f_h$ is their locally uniform limit.

If $f$ is equivariant and $h,k \in G$ then
\[f_{kh}(g_1,\ldots,g_p) = f_h(k^{-1}g_1,\ldots,k^{-1}g_p),\]
so if $(\g_n)_n$ is a sequence of type-I or type-II functions converging to $f_h$ then the functions $k^{-1}\cdot \g_n$ give a sequence of the same kind converging to $f_{kh}$.  Therefore type-I or type-II approximants for some $f_h$ can be used to give approximants for any other $f_{h'}$, so in this case the conclusion holds for every $h$. Finally, if $f$ is also regular at the identity, then we may choose the approximants $\g_n$ in the above construction to be locally constant around $(e,e,\ldots,e) \in G^{p+1}$, so that slicing each $\g_n$ at $e$ gives an approximant to $f_e$ which is locally constant around $(e,\ldots,e) \in G^p$.  Therefore $f_e$ is regular at the identity, and now the above equation implies also that $f_h$ is regular at $(h,\ldots,h)$. \qed

\begin{lem}
If $X$ is a locally compact and second countable metrizable space, $\mu$ is a Radon measure of full support on $X$ and $\cal{V}$ is an open cover of $X$, then there is a Borel partition $\cal{P}$ of $X$ such that
\begin{itemize}
\item $\cal{P}$ is locally finite;
\item each $P \in \cal{P}$ is contained in some member of $\cal{V}$;
\item and each $P \in \cal{P}$ satisfies $\mu(\partial P) = 0$.
\end{itemize}
\end{lem}

\noindent\textbf{Proof}\quad This construction rests on making careful use of a partition of unity; I doubt it is original, but have not found a suitable reference.

First, by local compactness we can express each $V \in \cal{V}$ as a union of precompact open subsets of $V$, and hence we may assume that every member of $\cal{V}$ is precompact.

By paracompactness we may choose a locally finite open refinement $\cal{U}$ of $\cal{V}$ and a partition of unity $(\rho_U)_U$ subordinate to $\cal{U}$.  Clearly it now suffices to prove the lemma with $\cal{U}$ in place of $\cal{V}$.  By second countability, $\cal{U}$ is countable.

Each member of $\cal{U}$ is precompact, and so by local finiteness there are values $\kappa_U \in (0,1)$ for each $U \in \cal{U}$ such that
\[\kappa_U < \frac{1}{|\{U' \in \cal{U}:\ U'\cap U \neq \emptyset\}|^2}.\]

If we now define $f:= \sum_U\kappa_U\rho_U:X\to \bbR$, then this is a strictly positive continuous function with the property that
\[f(x) \leq \sum_{U \in \cal{U}:\ U\ni x}\k_U < \frac{1}{|\{U \in\cal{U}:\ U \ni x\}|}\]
for all $x$.  This implies that for every $x \in X$ there is at least one $U \in \cal{U}$ for which $\rho_U(x) > f(x)$.  Therefore for any $s \in (0,1)$ the sets
\[Q_U^s := \{x \in X:\ \rho_U(x) > sf(x)\} \subseteq U\]
cover $X$, and this cover is also locally finite since each $Q_U^s$ is contained in its corresponding $U$. Moreover, for each fixed $U$ the boundaries
\[\partial Q_U^s \subseteq \{x \in X:\ \rho_U(x) = sf(x)\}, \quad s \in (0,1),\]
are pairwise disjoint, and so $\mu(\partial Q_U^s) = 0$ for Lebesgue-a.e. $s$. Since $\cal{U}$ is countable, it follows that there is some choice of $s \in (0,1)$ for which every $Q_U^s$ has boundary of measure zero.

Fix such an $s$ and let $Q_U := Q_U^s$.  Let $(Q_{U_i})_i$ be an enumeration of these sets, and for each $i$ let $P_i := Q_{U_i}\setminus \bigcup_{j< i} Q_{U_j}$.  Now $(P_i)_i$ is a locally finite Borel partition of $X$ having the desired properties. \qed

\begin{lem}[Equivariant continuation]\label{lem:equivt-cont}
In the setting of Lemma~\ref{lem:slice}, suppose that a function $f_0:G^p\to A$ is given which is almost type-I or, in case $G$ is an algebraic subgroup of $\rm{GL}_n(\bbR)$, almost type-II.  Then the same structure holds for the $G$-equivariant map $f:G^{p+1}\to A$ defined by
\[f(g_1,\ldots,g_p,g_{p+1}) := g_{p+1}\cdot \big(f_0(g_{p+1}^{-1}g_1,\ldots,g_{p+1}^{-1}g_p)\big).\]
If $f_0$ is regular at the identity then so is $f$.
\end{lem}

\noindent\textbf{Proof}\quad Let $(\eta_n)_n$ be a sequence of type-I (or type-II) functions converging locally uniformly to $f_0$ and define $G$-equivariant functions $\g_n:G^{p+1}\to A$ from each $\eta_n$ in the same way $f$ was defined from $f_0$.  Since the $G$-action on $A$ is continuous, these functions $\g_n$ converge locally uniformly to $f$, so it suffices to show that each $\g_n$ is itself an almost type-I (resp. almost type-II) function.  Note that $\g_n$ may not be \emph{exactly} type-I (resp. type-II), since the action of $g_{p+1}$ in its defining formula may give behaviour which is not locally constant.

Consider now a general l.c.s.c. group $G$ and a single type-I function $\eta:G^p\to A$. Since $\eta$ locally takes only finitely many values, every point $(h_1,\ldots,h_{p+1}) \in G^{p+1}$ has a precompact neighbourhood $V$ such that the function
\[\eta':(g_1,\ldots,g_{p+1})\mapsto \eta(g_{p+1}^{-1}g_1,\ldots,g_{p+1}^{-1}g_p)\]
takes only finitely many values on $V$.  Since the $G$-action on $A$ is continuous, for any $\eps > 0$ we may shrink $V$ further if necessary so that if $a_1$, \ldots, $a_\ell$ are these finitely many values then the sets
\[\{g_{p+1}\cdot a_i:\ (g_1,\ldots,g_{p+1}) \in V\}, \quad i=1,2,\ldots,\ell,\]
all have diameter less than $\eps$ in $A$, for some fixed choice of Polish metric on $A$.

Let $\cal{V}$ be a covering of $G^{p+1}$ by such neighbourhoods, and given this let $\cal{P}$ be the Borel partition obtained from $\cal{V}$ using the previous lemma.  Since any $P\in \cal{P}$ is contained in a member of $\cal{V}$, it admits a further partition $\cal{Q}_P$ into finitely many Borel subsets such that $\eta'$ is constant on each $Q \in \cal{Q}_P$ and
\[m_{G^{p+1}}(\partial Q) = 0 \quad \forall Q \in \cal{Q}_P.\]
Hence $\cal{Q} := \bigcup_P \cal{Q}_P$ is locally finite and consists of cells whose boundaries have measure zero, and by construction the map
\[\g(g_1,\ldots,g_{p+1}) := g_{p+1}\cdot \big(\eta'(g_1,\ldots,g_{p+1})\big)\]
is such that $\g(Q)$ has diameter less than $\eps$ in $A$ for every $Q \in \cal{Q}$.  Therefore if we let $\g'$ take a constant value from $\g(Q)$ on each of these sets $Q$, then $\g'$ is a type-I function that is $\eps$-uniformly close to $\g$, as required.

The case of an algebraic subgroup $G$ of $\rm{GL}_n(\bbR)$ and a type-II function $\eta$ is easier.  In that case we may always find a partition of $G^{p+1}$ which plays the r\^ole of the partition $\cal{P}$ above and consists of the intersections of $G$ with a partition of $M_{n\times n}(\bbR)\cong \bbR^{n^2}$ into dyadic cubes, which are manifestly semi-algebraic.  The rest of the argument is the same.

The last part of the conclusion is straightforward, since if $f_0$ is regular at the identity then in the above construction we can easily choose $\cal{P}$ and then $\cal{Q}$ such that the identity lies in the interior of its containing $\P$- and $\Q$-cells, so that the type-I or type-II approximants constructed above are locally constant around the identity. \qed

\vspace{7pt}

The heart of the inductive proof of Theorem B is the ability to lift functions of this type through quotient maps of target modules.

\begin{prop}[Lifting]\label{prop:lifting}
If $B\into A\onto A/B$ is an exact sequence of Polish Abelian groups, then any almost type-I function $f:G^p\to A/B$ which is regular at the identity has an almost type-I lift $G^p\to A$ which is regular at the identity.  If $G$ is algebraic then the same holds with `type-II' in place of `type-I'.
\end{prop}

\noindent\textbf{Proof}\quad Let $d$ be a translation-invariant Polish metric on $A$ and let $\ol{d}$ be the resulting quotient metric on $A/B$. Since $G^p$ is an l.c.s.c. group, it is also $\s$-compact, so we may choose an increasing sequence of compact subsets $K_n \subseteq G^p$ whose union is $G^p$.  Now let $(\g_n)_n$ be a sequence of type-I functions $G^p\to A/B$ such that
\[\ol{d}(\g_n(x),f(x)) < 2^{-n} \quad \forall x \in K_n\]
and with each $\g_n$ is locally constant around the identity.  Let $\P^0_n$ be the level-set partition of $\g_n$ and let $\P_n := \bigvee_{m\leq n}\P_m^0$, so each $\P_n$ is still a locally finite partition of $X$ with negligible boundary, each $\P_{n+1}$ is a refinement of $\P_n$, and for each $n$ the identity lies in the interior of its containing $\P_n$-cell.

Now one can recursively choose a sequence of lifts $\hat{\g}_n:G^p\to A$ of each $\g_n$ with the property that each $\hat{\g}_n$ is $\P_n$-measurable and
\[d(\hat{\g}_n(x),\hat{\g}_m(x)) \leq 2\ol{d}(\g_n(x),\g_m(x)) \quad \forall x.\]
To begin, let $\hat{\g}_1$ be any lift of $\g_1$ with the same level sets.  For the recursion, assume lifts $\hat{\g}_i$ have already been chosen for $i \leq n$.  For each $C \in \P_{n+1}$ we know that $\g_n$ and $\g_{n+1}$ are both constant on $C$.  If they are the same, then let $\hat{\g}_{n+1}$ take the same value as $\hat{\g}_n$ on $C$.  If they differ, then by the definition of the quotient metric we can choose $\hat{\g}_{n+1}(C)$ to be some element of $\g_{n+1}(C) + B$ that lies within distance $2\ol{d}(\g_n(C),\g_{n+1}(C))$ of $\hat{\g}_n(C)$ in $A$.

Each lift $\hat{\g}_n$ is still a type-I function and they satisfy the inequality
\[d(\hat{\g}_n(x),\hat{\g}_m(x)) < 4\cdot 2^{-m} \quad \hbox{whenever}\ n\geq m\ \hbox{and}\ x \in K_m,\]
so they form a locally uniformly Cauchy sequence.  Since $\hat{\g}_n$ is still $\P_n$-measurable, it is still locally constant at the identity. Letting $\hat{f}$ be the locally uniform limit of this sequence, it is an almost type-I function $G^p\to A$ which lifts $f$ and is regular at the identity.

The proof in case $G$ is algebraic and one wants almost type-II functions follows exactly the same steps. \qed

\begin{prop}\label{prop:type-I-and-II}
For any l.c.s.c. group $G$ and Polish $G$-module $A$, every cohomology class in $\rmH^p_\m(G,A)$ has a representative in the homogeneous bar resolution which is a $G$-equivariant almost type-I function $G^{p+1}\to A$ that is regular at the identity.  If, in addition, $G$ is an algebraic subgroup of some $\rm{GL}_n(\bbR)$, then this representative may be chosen to be almost type-II.
\end{prop}

\noindent\textbf{Proof}\quad We give the proof for general groups and almost type-I representatives, since the type-II case is almost identical now that Lemmas~\ref{lem:slice} and~\ref{lem:equivt-cont} have been proved.

This follows by an induction on degree using dimension-shifting.  When $p=0$ a cocycle is simply an element of $A^G$ regarded as a constant map $G\to A$, so is certainly of type-I or -II.  So now suppose the result is known for all degrees less than some $p\geq 1$ and that $\s:G^{p+1}\to A$ is a measurable cocycle.

Let $A' := \C(G,A)$. By dimension-shifting there is some $G$-equivariant $\psi:G^p\to A'$ such that $\s = d\psi$, where we identify $A$ with the subgroup of constant functions in $A'$.  Thus the map $\ol{\psi}:G^p\to A'/A$ obtained by quotienting is a cocycle, and so by the inductive hypothesis it is equal to $\ol{\phi} + d\ol{\kappa}$ for some almost type-I cocycle $\ol{\phi}:G^p\to A'/A$ that is regular at the identity and some $G$-equivariant measurable map $\ol{\kappa}:G^{p-1}\to A'/A$.

By Lemma~\ref{lem:slice} the slice
\[\ol{\phi}_0:(g_1,\ldots,g_{p-1})\mapsto \ol{\phi}(g_1,\ldots,g_{p-1},e)\]
is an almost type-I function on $G^{p-1}$ regular at the identity (if $p=1$ it is just a fixed element of $A'/A$).  Let $\phi_0:G^{p-1}\to A'$ be an almost type-I lift of it as promised by Proposition~\ref{prop:lifting}.  Lastly let $\phi:G^p\to A'$ be its equivariant continuation as in Lemma~\ref{lem:equivt-cont}, so this is also almost type-I and regular at the identity, and let $\kappa:G^{p-1}\to A'$ be any $G$-equivariant measurable lift of $\ol{\kappa}$ (such can always be found using the Measurable Selector Theorem).

Since $\psi$ is $G$-equivariant we know that
\[\psi = \phi + d\kappa + \a\]
for some equivariant $\a$ taking values in $A \leq A'$, so applying the differential gives
\[\s = d\phi + d\a.\]
It is easily seen from the alternating-sum formula for $d$ that $d\phi$ is still almost type-I and regular at the identity, and moreover the equation $d\phi = \s - d\a$ shows that it takes values in $A \leq A'$.  Any sequence $\eta_n$ of $A'$-valued type-I functions converging locally uniformly to $d\phi$ must therefore take values closer and closer to the subgroup $A$, and a small adjustment on each level set of each $\eta_n$ therefore gives a sequence of $A$-valued type-I functions converging locally uniformly to $d\phi$.  Thus $d\phi$ is an almost type-I $A$-valued representative for the cohomology class of $\s$ which is regular at the identity, and the induction continues. \qed

\vspace{7pt}

\noindent\textbf{Proof of Theorem B}\quad  If $\g_n:G^{p+1}\to A$ is a locally uniformly convergent sequence of type-I functions, and each $\g_n$ is locally constant on the full-measure open subset $U_n \subseteq G^{p+1}$, then $\lim_{n\to\infty}\g_n$ is still continuous on the full-measure G$_\delta$-set $\bigcap_nU_n$. \qed

\vspace{7pt}

\noindent\textbf{Proof of Theorem C}\quad If $A$ is discrete then a locally uniformly convergent sequence of type-I or type-II functions $\g_n:G^{p+1}\to A$ must eventually locally stabilize: that is, each point $x \in G^{p+1}$ has a neighbourhood $U$ such that all the restrictions $\g_n|_U$ are the same once $n$ is sufficiently large.  It follows that in this case the limits are still \emph{exactly} type-I or type-II.  Thus Proposition~\ref{prop:type-I-and-II} gives cocycle representatives that are type-I and, where applicable, type-II, and this is the content of Theorem C. \qed

\subsection*{The complex of locally continuous cochains}

The recent preprints~\cite{Fuc11a,Fuc11b,FucWoc11,WagWoc11} concern another variant of the bar resolution that can be used to compute a cohomology theory for topological groups.

Given a subset $U$ of $G$ and $p\geq 1$, let $\G_U^p$ denote the diagonal subset
\[\{(g_1,\ldots,g_{p+1}) \in G^{p+1}:\ g_i^{-1}g_j \in U\ \forall i\neq j\}.\]
Using these, one forms the complex of \textbf{locally continuous cochains}:
\begin{multline*}
\C_{\rm{lc}}^p(G,A):= \{\s\in \C(G^{p+1},A):\\ \exists\ \hbox{identity neighbourhood}\ U \subseteq G\ \rm{s.t.}\ \s|_{\G^p_U}\ \rm{continuous}\}.
\end{multline*}
Clearly this is a $G$-submodule of $\C(G^{p+1},A)$, and the alternating-sum differential $d$ satisfies $d(\C_\rm{lc}^p(G,A))\subseteq \C_\rm{lc}^{p+1}(G,A)$.  Cohomology groups $\rmH^\ast_\rm{lc}(G,A)$ may therefore be defined as the homology of the complex
\[0\to \C_\rm{lc}^0(G,A)^G\stackrel{d}{\to} \C_\rm{lc}^1(G,A)^G \stackrel{d}{\to} \C_\rm{lc}^2(G,A)^G \stackrel{d}{\to} \ldots. \]

Our definition of $\C^p_\rm{lc}(G,A)$ as a submodule of $\C(G^{p+1},A)$ implicitly restricts attention to measurable cochains, whereas Fuchssteiner, Wagemann and Wockel do not make this requirement.  However, some judicious measurable selection shows that this has no real effect on their results. Assuming that, the following theorem is a special case of results in~\cite{WagWoc11}.

\begin{thm}\label{thm:Seg-and-lc}
If $G$ is an l.c.s.c. topological group and $A$ is a topological $G$-module which is a k-space and locally contractible, then $\rmH_\Seg^\ast(G,A) \cong \rmH_\rm{lc}^\ast(G,A)$. \qed
\end{thm}

This is proved using a variant of Buchsbaum's criterion obtained in~\cite{WagWoc11} which gives a reduction to the case of a so-called `loop contractible' target module.  For that case, the works~\cite{Fuc11a,Fuc11b,FucWoc11} set up a spectral sequence relating $\rmH^\ast_\rm{lc}$ with the homology of the continuous  bar resolution (which correctly computes $\rmH^\ast_\Seg$ for a contractible module), which can be used to prove isomorphism of the continuous and locally-continuous theories in the necessary cases.

In the setting of l.c.s.c. groups and locally contractible Polish modules, the obvious inclusion $\l^p:\C_\rm{lc}^p(G,A) \subseteq \C(G^{p+1},A)$ immediately defines a connected sequence of comparison homomorphisms $\l^p_\ast:\rmH^p_\rm{lc}(G,A)\to \rmH^p_\m(G,A)$.  In view of Theorem~\ref{thm:Seg-and-lc}, another proof of Theorem A will result if one proves that each $\l^p_\ast$ is an isomorphism in case $A$ is discrete.

I do not know a quick proof of this, but at least the surjectivity of $\l^p_\ast$ follows at once from Theorem B.  That theorem tells us that any class in $\rmH^p_\m(G,A)$ has a representative $G^{p+1}\to A$ which is continuous at the identity, and so since $A$ is discrete it is actually locally constant on a neighbourhood of the identity.

By contrast, injectivity of $\l^p_\ast$ does not follow at once from Theorem B or C.  It requires one to prove that if $A$ is discrete, and if a locally continuous measurable cocycle $\s:G^{p+1}\to A$ is the boundary of a measurable cochain $\b:G^p\to A$, then $\b$ may also be chosen to be locally continuous.  However, I think this requires some result showing that classes in $\rmH^\ast_\rm{lc}$ also always have representative cocycles that have some useful addition structure, but this is already taking us closer to the proof of Theorem A in the following sections.

\section{Continuous dissections}\label{sec:diss-and-lay}

\subsection*{Continuous dissections}

A \textbf{dissection} of $I$ is a partition into finitely many intervals, all of them closed on the right and open on the left.

Henceforth $X$ will denote a metrizable topological space (the cases of interest will be $X = G^p$, $p \geq 1$).

\begin{dfn}[Continuous dissection; controlled partition]\label{dfn:cts-diss}
A \textbf{continuous dissection over $X$} is a family $\F$ of continuous functions $X\to [0,1]$ which contains the constant functions $0$ and $1$ and is \textbf{locally finite}, meaning that every $x \in X$ has a neighbourhood $U$ such that the set $\{\xi|_U:\ \xi \in\cal{F}\}$ is finite.

If $\F$ is a continuous dissection then an \textbf{$\F$-wedge} is a subset of $X\times I$ of the form
\[\{(x,t):\ \xi_1(x) < t \leq \xi_2(x)\}\]
for some $\xi_1,\xi_2 \in \F$.  A partition $\P$ of $X\times I$ is \textbf{controlled} by $\F$ if each of its cells is a union of $\F$-wedges.
\end{dfn}

\begin{figure}
\centering
\includegraphics[width=0.9\textwidth]{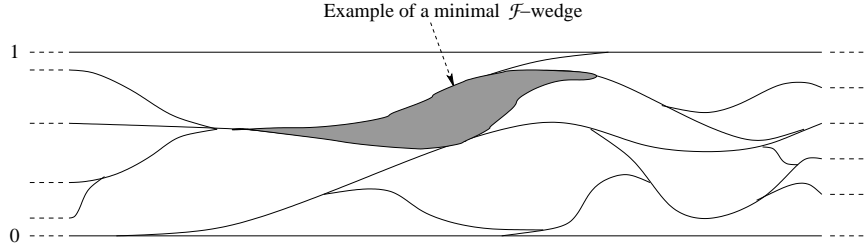}
\caption{Part of a continuous dissection over $\bbR$.  The $\F$-wedge shown includes its upper boundary, but not its lower.}\label{fig:cts-diss}
\end{figure}

Figure~\ref{fig:cts-diss} sketches an example of a continuous dissection $\F$ over $\bbR$, and highlights one of the resulting $\F$-wedges.

By the local finiteness of $\F$ and the continuity of its members, we may think of $\{\xi(x):\ \xi \in \F\}$ as specifying the end-points of a dissection of $I$ that varies continuously with $x$.  This motivates the terminology.

Clearly the union of any finite family of continuous dissections is still a continuous dissection.

If $\zeta,\xi:X\to \bbR$ are continuous functions then $\zeta\vee \xi$ and $\zeta\wedge \xi$ will denote their pointwise maximum and pointwise minimum respectively.

\begin{lem}
If $\F$ is a continuous dissection over $X$ then so is the family $\ol{\F}$ consisting of all functions obtained from members of $\F$ by repeated applications of $\wedge$, $\vee$ and pointwise limits of convergent directed families.
\end{lem}

\noindent\textbf{Proof}\quad Any maximum or minimum of continuous functions is still continuous, and if $U \subseteq X$ is open and such that $\{\xi|_U:\ \xi \in \F\}$ is finite, then
\[\{\zeta|_U:\ \zeta \in \ol{\F}\} = \ol{\{\xi|_U:\ \xi \in \F\}}\]
is still finite. \qed

\begin{dfn}[Lattice-completeness]
The continuous dissection $\ol{\F}$ constructed from $\F$ as above is the \textbf{lattice-hull} of $\F$, and $\F$ itself is \textbf{lattice-complete} (`\textbf{l-complete}') if $\F = \ol{\F}$.
\end{dfn}

Observe that if $C$ is an $\F$-wedge then $(X\times I)\setminus C$ is either empty, an $\F$-wedge or a union of two $\F$-wedges.  This easily implies the following.

\begin{lem}
If $\F$ is l-complete then any nonempty intersection of $\F$-wedges is an $\F$-wedge, and hence each point of $X\times I$ lies in a unique minimal $\F$-wedge. The minimal $\F$-wedges define a locally finite partition of $X\times I$. \qed
\end{lem}

Continuous dissections behave well under pulling back.

\begin{lem}[Pulling back continuous dissections]\label{lem:CD-pullback}
If $\phi:X\to Y$ is a continuous map between metrizable spaces and $\F$ is a continuous dissection over $Y$, then the family 
\[\phi^\ast\F := \{\xi\circ \phi:\ \xi \in \F\}\]
is a continuous dissection over $X$.
\end{lem}

\noindent\textbf{Proof}\quad Continuity of each $\xi\circ\phi$ is immediate, and the local finiteness of $\phi^\ast\F$ follows because for any $x \in X$ there is a neighbourhood $U$ of $\phi(x)$ on which $\F$ restricts to a finite family, and now by continuity $\phi^{-1}(U)$ is a neighbourhood of $x$ on which $\phi^\ast\F$ restricts to a finite family. \qed

\vspace{7pt}

Much of the versatility of continuous dissections derives from the following construction (and its relative in Lemma~\ref{lem:ctsdiss-loc-to-glob-approx} below).

\begin{lem}\label{lem:wedges-in-cover}
If $\U$ is an open cover of $X$ then there is a continuous dissection $\F$ over $X$ such that every minimal $\F$-wedge is contained in $U\times I$ for some $U \in \U$.
\end{lem}

\noindent\textbf{Proof}\quad By paracompactness we may assume that $\U$ is locally finite and choose a subordinate partition of unity $(\rho_U)_U$. Now let $\F$ be the class of all functions of the form
\[\tau_{U_1,\ldots,U_m} := \rho_{U_1} + \ldots + \rho_{U_m}\]
for some $U_1$, \ldots, $U_m \in \U$.

These are continuous and $[0,1]$-valued, and $\F$ is clearly locally finite, so it is a continuous dissection.  

Suppose that $(x,t) \in X\times I$. Then since $(\rho_U)_U$ is a partition of unity, there are some distinct $U_1$, $U_2$, \ldots, $U_m \in \cal{U}$ such that
\[\rho_{U_1}(x) + \cdots + \rho_{U_{m-1}}(x) < t \leq \rho_{U_1}(x) + \cdots + \rho_{U_m}(x).\]
Letting $\tau_1(x)$ and $\tau_2(x)$ denote the members of $\F$ appearing on the left- and right-hand sides here, we have shown that
\[(x,t) \in \{(x',t'):\ \tau_1(x') < t' \leq \tau_2(x')\} \subseteq \{(x',t'):\ \rho_{U_m}(x') > 0\} \subseteq U_m\times I,\]
so $(x,t)$ is contained in an $\F$-wedge which is itself contained in the lift of a member of $\cal{U}$.  Since $(x,t)$ was arbitrary, all minimal $\F$-wedges must have this property, as required. \qed

\subsection*{Product spaces and ascending tuples}

In general we will need to handle functions defined on spaces of the form
\[X_1\times \cdots \times X_p \times I^p\]
for some metrizable spaces $X_1$, \ldots, $X_p$, $p \geq 1$.  This will involve working with whole $p$-tuples of continuous dissections, in which the $i^{\rm{th}}$ continuous dissection applies to the $i^{\rm{th}}$ coordinate in $I^p$ for $i=1,2,\ldots,p$.  Moreover, it will be crucial that these tuples of continuous dissections respect the product structure of $X_1\times \cdots \times X_p$ in the following very particular way.

\begin{dfn}[Ascending tuples]
If $X_1$, \ldots, $X_p$ is a tuple of metrizable spaces, then a tuple of continuous dissections $(\F_1,\ldots,\F_p)$ is \textbf{ascending} over $X_1$, \ldots, $X_p$ if
\begin{quote}
$\F_1$ is a continuous dissection over $X_1$,

$\F_2$ is a continuous dissection over $X_1\times X_2$,

\quad\quad\quad\quad $\vdots$

$\F_p$ is a continuous dissection over $X_1 \times \cdots \times X_p$.
\end{quote}
\end{dfn}

In the sequel, when a tuple of spaces $X_1$, \ldots, $X_p$ is understood, we will usually abbreviate
\[X_{\leq i} := X_1\times \cdots \times X_i \quad \hbox{for}\ i=1,2,\ldots,p.\]
Occasionally we will have need for the coordinate projections $X_{\leq j}\to X_{\leq i}$ for $i < j$.  We denote these by $\pi_{\leq i}$, since the dependence on $j$ should always be clear.

\begin{dfn}[Multiwedges and control]\label{dfn:multi-control}
If $(\F_1, \ldots, \F_p)$ is an ascending tuple of continuous dissections over $X_1$, \ldots, $X_p$, then an \textbf{$(\F_1,\ldots,\F_p)$-multiwedge} is a subset of $X_{\leq p}\times I^p$ of the form
\[\{(x_1,\ldots,x_p,t_1,\ldots,t_p):\ (x_1,\ldots,x_i,t_i) \in C_i\ \forall i=1,2,\ldots,p\},\]
where $C_i$ is an $\F_i$-wedge for each $i$.  This multiwedge will sometimes be written as the fibred product
\[C_1 \times_{X_{\leq p}} C_2 \times_{X_{\leq p}} \cdots \times_{X_{\leq p}} C_p\]
(this is slightly abusive, since formally the wedges $C_i$ are defined over the different spaces $X_{\leq i}$, but no confusion will arise).

A partition $\P$ of $X_{\leq p}\times I^p$ is \textbf{controlled by $(\F_1,\ldots,\F_p)$} if every cell of $\P$ is a union of $(\F_1,\ldots,\F_p)$-multiwedges.
\end{dfn}

\begin{lem}
If each $\F_i$ is l-complete and a given $(\F_1,\ldots,\F_p)$-multiwedge is minimal under inclusion, then it can be expressed as the fibred product of minimal $\F_i$-wedges.
\end{lem}

\noindent\textbf{Proof}\quad If $(x_1,\ldots,x_p,t_1,\ldots,t_p) \in X_{\leq p}\times I^p$, then an easy check shows that the minimal $(\F_1,\ldots,\F_p)$-multiwedge containing it must be equal to
\[C_1 \times_{X_{\leq p}}\cdots \times_{X_{\leq p}}C_p,\]
where each $C_i$ is the minimal $\F_i$-wedge containing $(x_1,\ldots,x_i,t_i)$. \qed

\vspace{7pt}

Henceforth we will always assume that our continuous dissections are l-complete.

Ascending tuples also enjoy an analog of Lemma~\ref{lem:CD-pullback} in terms of the following class of maps.

\begin{dfn}[Ascending maps]
Suppose that $X_i$ and $Y_i$ for $i = 1,2,\ldots,p$ are metrizable spaces.  Then an \textbf{ascending tuple of maps} from $X_1$, \ldots, $X_p$ to $Y_1$, \ldots, $Y_p$ is a tuple of continuous maps
\begin{quote}
$\phi_1:X_1\to Y_1$,

$\phi_2:X_{\leq 2}\to Y_2$,

\quad\quad\quad\quad $\vdots$

$\phi_p:X_{\leq p}\to Y_p$.
\end{quote}

Given these, we will define further maps $\phi_{\leq i}:X_{\leq i}\to Y_{\leq i}$ for $i=1,2,\ldots,p$ by
\[\phi_{\leq i}:X_{\leq i}\to Y_{\leq i}:
(x_1,\ldots,x_i)\mapsto \big(\phi_1(x_1),\ldots,\phi_i(x_1,\ldots,x_i)\big).\]
\end{dfn}

The following extension of Lemma~\ref{lem:CD-pullback} is immediate.

\begin{lem}
If $\F_1$, \ldots, $\F_p$ is an ascending tuple of continuous dissections over $Y_1$, \ldots, $Y_p$, and $\phi_i:X_{\leq i}\to Y_i$ is an ascending tuple of maps, then the tuple of continuous dissections
\[\phi_1^\ast\F_1,\ \phi_{\leq 2}^\ast\F_2,\ \ldots,\ \phi_{\leq p}^\ast\F_p\]
is ascending over $X_1$, \ldots, $X_p$. \qed
\end{lem}

The last result of this section is a technical extension of Lemma~\ref{lem:wedges-in-cover} that will be crucial later.

\begin{lem}\label{lem:ctsdiss-loc-to-glob-approx}
Suppose that $\cal{U}$ is an open cover of $X_1$ and that to every $U \in \cal{U}$ there is associated an ascending tuple $\F_{U,1},\ldots,\F_{U,p}$ of continuous dissections over $X_1,\ldots,X_p$. Then there is another ascending tuple $\F_1,\ldots,\F_p$ with the following property: for every minimal $(\F_1,\ldots,\F_p)$-multiwedge $C$ there is some $U \in \cal{U}$ such that
\begin{itemize}
\item $C \subseteq U\times X_2\times \cdots \times X_p \times I^p$, and
\item $C$ is contained in some $(\F_{U,1},\ldots,\F_{U,p})$-multiwedge.
\end{itemize}
\end{lem}

\noindent\textbf{Remark}\quad As the proof will show, it is essential that the open sets here $U$ depend only on the coordinate in $X_1$. \fin

\vspace{7pt}

\noindent\textbf{Proof}\quad By paracompactness we may assume that $\cal{U}$ is locally finite.  Having done so, another quick appeal to paracompactness gives a further locally finite refinement $\cal{V}$ of $\cal{U}$ such that for each $V \in \cal{V}$ the collection
\[\cal{U}_V := \{U \in \cal{U}:\ V\cap U \neq \emptyset\}\]
is finite.

Now for each $V \in \cal{V}$ we choose a $U_V \in \cal{U}$ that contains it, and set
\[(\F_{V,1},\ldots,\F_{V,p}) := (\F_{U_V,1},\ldots,\F_{U_V,p}).\]
By local finiteness, we may let $(\rho_V)_V$ be a partition of unity subordinate to $\cal{V}$, and now as in Lemma~\ref{lem:wedges-in-cover} let $\cal{G}_1$ be the continuous dissection over $X_1$ given by the lattice-hull of the functions $0$, $1$ and
\[\tau_{V_1,\ldots,V_m} := \rho_{V_1} + \ldots + \rho_{V_m} \quad \hbox{for}\ V_1,\ldots,V_m \in \cal{V}.\]
Just as in Lemma~\ref{lem:wedges-in-cover}, it follows that every minimal $\cal{G}_1$-wedge is contained in a set of the form $V\times I$ for some $V \in \cal{V}$.  Also, for $i=2,\ldots,p$ let $\cal{G}_i$ be the pullback of $\cal{G}_1$ through the coordinate projection $\pi_1:X_{\leq i}\to X_1$, and let $\tau^{(i)}_{V_1,\ldots,V_m} := \tau_{V_1,\ldots,V_m}\circ\pi_1$.

Finally, for $i=1,2,\ldots,p$ we define
\[\F_i = \ol{\{0,1\}\cup\bigcup_{V_1,\ldots,V_m \in \cal{V},\ U \in \cal{U}\ \rm{s.t.}\ V_m \cap U\neq\emptyset}\big\{\tau^{(i)}_{V_1,\ldots,V_{m-1}}\vee(\xi\wedge \tau^{(i)}_{V_1,\ldots,V_m}):\ \xi \in \F_{U,i}\big\}}\]
(one checks easily that this is locally finite).

This is an ascending tuple over $X_1$, \ldots, $X_p$.  We will show that it has the desired two properties.  Suppose that
\[C = C_1 \times_{X_{\leq p}} \cdots \times_{X_{\leq p}} C_p\]
is a minimal $(\F_1,\ldots,\F_p)$-multiwedge.  We may write this representation so that each $C_i$ is a minimal $\F_i$-wedge, and so since $\cal{G}_i \subseteq \F_i$, each $C_i$ must lie in some $\cal{G}_i$-wedge of the form
\begin{multline*}
D_i := \{(x_1,\ldots,x_i,t_i):\ \tau_{V^i_1,\ldots,V^i_{m_i - 1}}(x_1) < t_i \leq \tau_{V^i_1,\ldots,V^i_{m_i}}(x_1)\}\\ \subseteq V^i_{m_i}\times X_2\times \cdots\times X_i\times I,\end{multline*}
implying that
\[C_i \subseteq V^i_{m_i}\times X_2\times \cdots \times X_i\times I.\]

Suppose $(x_1,\ldots,x_p,t_1,\ldots,t_p) \in C$.  Then for each $i$ we have
\[(x_1,\ldots,x_i,t_i) \in C_i,\]
which requires in particular that $x_1 \in V^i_{m_i}$ for $i=1,2,\ldots,p$.  This implies that if $U = U_{V^1_{m_1}}$, then $U$ still has nonempty intersection with $V^i_{m_i}$ for all $i=2,\ldots,p$. 

Now on the one hand we have
\[C \subseteq V^1_{m_1}\times X_2\times \cdots \times X_p\times I^p \subseteq U \times X_2\times \cdots \times X_p\times I^p,\]
which proves the first property. On the other hand, within the $\cal{G}_i$-wedge $D_i$ introduced above, the partition into minimal $\F_i$-wedges is a refinement of the partition into 
minimal $\F_{U',i}$-wedges for any $U' \in \cal{U}$ that intersects $V^i_{m_i}$.  Our choice of $U$ above is one such member of $\cal{U}$, so our minimal $(\F_1,\ldots,\F_p)$-multiwedge $C$ is contained in some minimal $(\F_{U,1},\ldots,\F_{U,p})$-wedge, as required for the second property. \qed

\section{Semi-layered functions}\label{sec:sl}

\subsection*{Layered and semi-layered functions}

Now suppose that $A$ is a Hausdorff topological group.  In the coming application to cohomology, our interest will be in $A$-valued functions on Cartesian powers $G^p$ of a l.c.s.c. group $G$.  In this setting, it will be important that we work with a class of functions that respects the order of the coordinate factors in $G^p$.

More generally, suppose again that $X_1$, \ldots, $X_p$ are metrizable spaces.  Let $\b:X_{\leq p}\times I^p\to X_{\leq p}$ be the obvious coordinate projection between these spaces.  Later we will focus on the case $X_1 = \cdots = X_p = G$, but the order of the coordinates will still be important. The class of functions we need is the following.

\begin{dfn}[Layered and semi-layered functions]\label{dfn:layered} For a given tuple of spaces $X_1$, \ldots, $X_p$, a function
\[\g:X_{\leq p}\times I^p\to A\]
is \textbf{layered} if there is an ascending tuple of l-complete continuous dissections $\F_1$, \ldots, $\F_p$ over $X_1$, \ldots, $X_p$ such that $\g$ is constant on every minimal $(\F_1,\ldots,\F_p)$-multiwedge.  In this case we write that $\g$ itself is \textbf{controlled} by $(\F_1,\ldots,\F_p)$.

Similarly, a function
\[f:X_{\leq p}\times I^p\to A\]
is \textbf{semi-layered} if there is such an ascending tuple $\F_1$, \ldots, $\F_p$ such that for every minimal $(\F_1,\ldots,\F_p)$-multiwedge $C$ there is a continuous function $f_C:\ol{\b(C)} \to A$ such that
\[f|_C = f_C\circ \b|_C\]
(so, in particular, $f|_C(x_1,\ldots,x_p,t_1,\ldots,t_p)$ does not depend on $(t_1,\ldots,t_p)$ when $(x_1,\ldots,x_p,t_1,\ldots,t_p)$ is assumed to lie in $C$). In this case we write that $f$ is \textbf{semi-controlled} by $(\F_1,\ldots,\F_p)$.
\end{dfn}

Note that the definitions of layered and semi-layered functions make implicit reference to the structure of $X_{\leq p}$ as a product of the spaces $X_1$, \ldots, $X_p$.

\vspace{7pt}

\noindent\textbf{Example}\quad If $p = 1$, a function $f:X_1\times I\to A$ is layered if it is constant on each minimal $\F_1$-wedge (recall the sketch in Figure~\ref{fig:cts-diss}). It is semi-layered if for each minimal $\F_1$-wedge $C$, $f|_C$ is lifted from some continuous function on $\ol{\b(C)}$. \fin

\vspace{7pt}

\noindent\textbf{Example}\quad Suppose that $f:X_{\leq p}\times I^p\to A$ is layered and controlled by a tuple $(\F_1,\ldots,\F_{p-1},\F)$ in which $\F_i$ is the trivial continuous dissection $\{0,1\}$ for all $i\leq p-1$.  Then $f$ has no dependence on the first $p-1$ copies of $I$, and it may be regarded as a layered function $X' \times I\to A$ in the case $p = 1$ controlled by $\F$, where the product structure of $X' := X_{\leq p}$ is forgotten.  This simple observation will be useful shortly. \fin

\vspace{7pt}

It is easy to show that any layered function $f$ is also semi-layered, but in this case semi-control by a tuple $(\F_1,\ldots,\F_p)$ does not imply control by the same tuple.  Specifically, if $C$ is an $(\F_1,\ldots,\F_p)$-multiwedge for which $\ol{\b(C)} \subseteq X_{\leq p}$ has more than one connected component, then $f$ could take different values on the lifts of those components and still be lifted from a continuous function on $\ol{\b(C)}$.

The following is immediate.

\begin{lem}
If $\cal{G}_i \supseteq \F_i$ are continuous dissections as above for each $i$ and $\g:X_{\leq p}\times I^p\to A$ is a layered (resp. semi-layered) function controlled (resp. semi-controlled) by $(\F_1,\ldots,\F_p)$, then it is also controlled (resp. semi-controlled) by $(\cal{G}_1,\ldots,\cal{G}_p)$. \qed
\end{lem}

\begin{lem}\label{lem:sums-still-s-l}
The sum of two (semi-)layered functions is (semi-)layered.
\end{lem}

\noindent\textbf{Proof}\quad If $f_1,f_2:X_{\leq p}\times I^p\to A$ are (semi-)layered and are respectively (semi-)controlled by $(\F^1_1,\ldots,\F^1_p)$ and $(\F^2_1,\ldots,\F^2_p)$, then $f_1 + f_2$ is (semi-)controlled by $(\ol{\F^1_1\cup \F^2_1},\ldots,\ol{\F^1_p\cup \F^2_p})$. \qed

\vspace{7pt}

Layered functions also exhibit good behaviour under pulling back.  The correct formulation of this behaviour is a little delicate.

\begin{lem}[Pulling back and slicing]\label{lem:pre-pullback}
Suppose that $\phi_i:X_{\leq i}\to Y_i$ is an ascending tuple of maps between metrizable spaces and that $\g:Y_{\leq p}\times I^p\to A$ is a layered function controlled by $(\F_1,\ldots,\F_p)$.  Abbreviate $\phi_{\leq p} =: \phi$. Then the pullback $\phi^\ast\g := \g(\phi(\cdot),\cdot)$ is a layered function on $X_{\leq p}\times I^p$, controlled by the tuple $(\phi_1^\ast\F_1,\ldots,\phi_{\leq p}^\ast\F_p)$.
The analogous assertion holds for semi-layered functions and semi-control.
\end{lem}

\noindent\textbf{Proof}\quad Both conclusions follow from the behaviour of the pulled-back continuous dissections. Suppose that
\[C = \{(y_1,\ldots,y_p,t_1,\ldots,t_p):\ (y_1,\ldots,y_i,t_i) \in C_i\ \forall i\leq p\}\]
is a minimal $(\F_1,\ldots,\F_p)$-multiwedge with each $C_i$ being a minimal $\F_i$-wedge.  Then the pullback of this set under $\phi\times \id_{I^p}$ is a multiwedge for the tuple $(\phi_1^\ast\F_1,\ldots,\phi_{\leq p}^\ast\F_p)$.  Hence if $\g:Y_{\leq p}\times I^p\to A$ is layered and controlled by the $\F_i$, then its pullback $\phi^\ast\g$ is controlled by these pullbacks $\phi_i^\ast\F_i$.  On the other hand, if $f$ is semi-layered and semi-controlled by these $\F_i$, then for each minimal $(\F_1,\ldots,\F_p)$-multiwedge $C$ there is a continuous function $f_C:\ol{\b(C)}\to A$ such that $f|_C = f_C\circ \b|_C$.  This now implies
\[\phi^\ast f|_{(\phi\times \id_{I^p})^{-1}(C)} = (f_C\circ \phi)\circ \b|_{(\phi\times \id_{I^p})^{-1}(C)},\]
where $f_C\circ \phi$ is a continuous function defined on the set
\[\phi^{-1}(\ol{\b(C)}) \supseteq \ol{\b((\phi\times \id_{I^p})^{-1}(C))}.\]
So the conditions of the second part of Definition~\ref{dfn:layered} are still satisfied. \qed

\vspace{7pt}

We next present the key analytic result that will give us some control over the possible discontinuities of cocycles, by applying it during an induction by dimension-shifting.  Its proof illustrates the use of Lemma~\ref{lem:wedges-in-cover}.

\begin{prop}[Lifting semi-layered functions]\label{prop:lift-s-l}
Suppose that $B \into A \onto A/B$ is an exact sequence of Hausdorff topological Abelian groups that admits a local continuous cross-section.  Then any semi-layered function $f:X_{\leq p}\times I^p \to A/B$ has a semi-layered lift $X_{\leq p}\times I^p\to A$.
\end{prop}

\noindent\textbf{Proof}\quad Let $f:X_{\leq p}\times I^p\to A/B$ be a semi-layered function, and let $\cal{P}$ be the partition of $X_{\leq p}\times I^p$ into minimal $(\F_1,\ldots,\F_p)$-multiwedges.  Let $\b:X_{\leq p}\times I^p\to X_{\leq p}$ be the coordinate projection.

Since each $\F_i$ is locally finite, any $x \in X_{\leq p}$ can intersect only finitely many of the closures $\ol{\b(C)}$ with $C \in \P$.  Having fixed such a point $x$, let $C_1$, \ldots, $C_\ell$ be these members of $\P$, and for each $i \leq \ell$ let $f_i:\ol{\b(C_i)}\to A/B$ be a continuous map such that $f|_{C_i} = f_i\circ\b|_{C_i}$.

Since $A \onto A/B$ admits continuous local sections, for each $i \leq \ell$ we can choose a neighbourhood $V_i$ of $f_i(x)$ such that $\ol{V_i}$ admits such a section. For each $i$, $f_i^{-1}(V_i)$ is a relatively open subset of $\ol{\b(C_i)}$ containing $x$, so we may find a neighbourhood $U_{x,i}$ of $x$ such that $U_{x,i}\cap \ol{\b(C_i)} \subseteq f^{-1}(V_i)$, and now $U_x := \bigcap_{i\leq \ell}U_{x,i}$ is still a neighbourhood of $x$.

The neighbourhoods $U_x$ obtained this way comprise an open cover of $X_{\leq p}$, so Lemma~\ref{lem:wedges-in-cover} gives an l-complete continuous dissection $\cal{G}_0$ over $X_{\leq p}$ such that every minimal $\cal{G}_0$-wedge is contained in $U_x \times I$ for some $x$.  Letting $\cal{G} := \ol{\cal{G}_0 \cup \cal{F}_p}$, it follows that
\begin{itemize}
\item any minimal $(\F_1,\ldots,\F_{p-1},\cal{G})$-multiwedge $D$ is both contained in some minimal $(\F_1,\ldots,\F_p)$-multiwedge $C$, and also in $\b^{-1}(U_x)$ for some $x$, and
\item if $D\subseteq C$ are as above and $f|_C = f_C\circ \b|_C$ with $f_C:\ol{\b(C)}\to A/B$ continuous, then $f_C(\ol{\b(D)})$ is contained in an open subset of $A/B$ that admits a continuous section to $A$.
\end{itemize}

Let $\Phi_D:f_C(\ol{\b(D)})\to A$ be such a continuous section for each $D$, and define the function $F:X_{\leq p}\times I^p\to A$ by
\[F|_D = \Phi_D\circ (f|_D) \quad \forall\ \hbox{minimal $(\F_1,\ldots,\F_{p-1},\cal{G})$-multiwedge}\ D.\]
For each minimal $(\F_1,\ldots,\F_{p-1},\cal{G})$-multiwedge $D$ the restriction $F|_D$ is given by $(\Phi_D\circ f_C)\circ \b|_D$, where $C$ is the minimal $(\F_1,\ldots,\F_p)$-multiwedge containing $D$, and the function $\Phi_D\circ f_C:\ol{\b(D)}\to A$ is continuous.  Therefore $F$ is a semi-layered lift of $f$, semi-controlled by $(\F_1,\ldots,\F_{p-1},\cal{G})$.  \qed

\subsection*{Semi-layered functions and Segal's soft modules}

Definition~\ref{dfn:layered} is motivated by the need to define a `concrete' class of functions on $G^p\times I^p$ that lie within the modules appearing in Segal's resolution~(\ref{eq:Seg-cplx}).  The following lemma tells us that semi-layered functions form such a class. In practice, it will be used to show that the `semi-layered' cohomology theory is effaceable in Segal's category.

\begin{lem}\label{lem:pre-sl-are-Seg}
If $f:X\times I\to A$ is semi-layered then setting
\[F(x)(\cdot) := f(x,\cdot)\]
defines a continuous function $X\to EA$ (that is, an element of $\C_\cts(X,EA)$).
\end{lem}

\noindent\textbf{Proof}\quad Suppose that $f$ is semi-controlled by the l-complete continuous dissection $\F$ and let $\P$ be the partition of $X\times I$ into minimal $\F$-wedges. Each $x \in X$ has a neighbourhood $U$ such that $\{\xi|_U:\ \xi \in \F\}$ is finite, so we may enumerate this set of restricted functions as $\xi_1$, \ldots, $\xi_m$.  Also, $x$ can lie in the closure of $\ol{\b(C)}$ for only finitely many $C \in \P$, say $C_1$, \ldots, $C_r$, and for each of these there is a continuous function $f_j:\ol{\b(C_j)}\to A$ such that $f|_{C_j} = f_j\circ \b|_{C_j}$.

By continuity, given $\eps > 0$ and an identity neighbourhood $V$ in $A$, we may now shrink $U$ further if necessary so that
\begin{itemize}
\item there are values $t_1$, \ldots, $t_m \in [0,1]$ such that $|\xi_i(y) - t_i| < \eps$ for each $i \leq m$ and $y \in U$, and
\item $f_j(y) \in f_j(x) + V$ for all $j \leq r$ and $y \in U$.
\end{itemize}

These conditions imply that $f(y,\cdot)$ lies within a small neighbourhood of $f(x,\cdot)$ in $EA$ for all $y \in U$; since $\eps$ and $V$ were arbitrary, this completes the proof. \qed

\vspace{7pt}

\begin{cor}\label{cor:fn-vald}
If $f:X_{\leq p}\times I^p\to A$ is semi-layered then the function $F:X_{\leq p-1}\times I^{p-1}\to A^{X_p\times I}$ defined by
\[F(x_1,\ldots,x_{p-1},t_1,\ldots,t_{p-1})(\cdot,\cdot) := f(x_1,\ldots,x_{p-1},\cdot,t_1,\ldots,t_{p-1},\cdot)\]
takes values in $\C_\cts(X_p,EA)$.
\end{cor}

\noindent\textbf{Proof}\quad Let $\b:X_{\leq p}\times I^p\to X_{\leq p}$ and $\b_p:X_p\times I\to X_p$ be the coordinate projections.

In order to apply the previous lemma, we need to show that for every $x_1$, \ldots, $x_{p-1}$, $t_1$, \ldots, $t_{p-1}$ the function
\[f(x_1,\ldots,x_{p-1},\cdot,t_1,\ldots,t_{p-1},\cdot)\]
is semi-layered.  To see this, suppose that $f$ is semi-controlled by $(\F_1,\ldots,\F_p)$, and fix $(x^-,t^-) := (x_1,\ldots,x_{p-1},t_1,\ldots,t_{p-1})$.  Define
\[\cal{G} := \{\xi(x^-,\cdot):\ \xi \in \F_p\},\]
so that any minimal $\cal{G}$-wedge $D$ is of the form
\[\{(x,t):\ \xi_1(x^-,x) < t \leq \xi_2(x^-,x)\} \quad \hbox{for some}\ \xi_1,\xi_2 \in \F_p.\]
This can be identified with $(\{x^-\}\times X_p \times I)\cap C_{p,D}$ for some choice of $\F_p$-wedge $C_{p,D}$, which may also be assumed to be minimal.

Let $C_i$ for $i=1,\ldots,p-1$ be minimal $\F_i$-wedges such that
\[(x^-,t^-) \in C^- := C_1 \times_{X_{\leq p-1}} \cdots \times_{X_{\leq p-1}} C_{p-1}.\]
Then for any minimal $\cal{G}$-wedge $D$ one has
\[\{(x^-,t^-)\}\times D = (\{(x^-,t^-)\}\times (X_p\times I)\big) \cap C_D\]
where
\[C_D = C_1 \times_{X_{\leq p}} \cdots \times_{X_{\leq p}} C_{p-1}\times_{X_{\leq p}} C_{p,D},\]
which is a minimal $(\F_1,\ldots,\F_p)$-multiwedge.  Therefore Definition~\ref{dfn:layered} gives a continuous function $f_{C_D}:\ol{\b(C_D)}\to A$ such that $f|_{C_D} = f_{C_D}\circ\b|_{C_D}$, and in particular
\[f(x^-,x,t^-,t) = f_{C_D}(x^-,x) \quad \forall (x,t) \in D.\]

This may now be re-written as
\[F(x^-,t^-)(\cdot,\cdot)|_D = f_{C_D}(x^-,\cdot) \circ \b_p|_D.\]
To finish the proof, observe that if $x \in \ol{\b_p(D)}$, then for any neighbourhood $U$ of $x$ in $X_p$ there is some $x' \in U\cap \b_p(D)$, and hence there is also some $t' \in I$ such that
\[(x^-,x',t') \in C_{p,D} \quad \Longrightarrow \quad (x^-,x',t^-,t') \in C_D.\]
Therefore $(x^-,x') \in \b(C_D)$, and since $x'$ was arbitrarily close to $x$ it follows that $(x^-,x) \in \ol{\b(C_D)}$.  Hence $\{x^-\}\times \ol{\b_p(D)} \subseteq \ol{\b(C)}$, and so we may define
\[F_D(x) := f_{C_D}(x^-,x) \quad \hbox{for}\ x \in \ol{\b_p(D)}.\]
This gives a continuous function $F_D:\ol{\b_p(D)}\to A$ such that $F(x^-,t^-)|_D = F_D\circ\b_p|_D$, and so proves that $F(x^-,t^-)(\cdot,\cdot)$ is semi-layered for each fixed $(x^-,t^-)$, as required. \qed

\vspace{7pt}

We will also need the following enhancement to the above corollary.

\begin{prop}\label{prop:fn-vald-semi-layered}
If $f$ and $F$ are as in the preceding corollary and $f$ is semi-controlled by $(\F_1,\ldots,\F_p)$, then $F$ is a semi-layered as a function $X_{\leq p-1}\times I^{p-1}\to \C_\cts(X_p,EA)$ and is semi-controlled by $(\F_1,\ldots,\F_{p-1})$.
\end{prop}

The proof of this will use two auxiliary lemmas.

\begin{lem}\label{lem:pre-fn-vald-semi-layered}
Let $X$ and $Y$ be metrizable spaces and $A$ a Hausdorff topological group, and suppose that $f:(X\times Y)\times I\to A$ is a semi-layered function if we ignore the product structure of $X\times Y$, semi-controlled by an l-complete continuous dissection $\F$ over $X\times Y$.  Then the map
\[F:x\mapsto f(x,\cdot,\cdot)\]
takes values in $\C_\cts(Y,EA)$ and is continuous for the Segal topology on that module.
\end{lem}

\noindent\textbf{Proof}\quad That $F$ takes values in $\C_\cts(Y,EA)$ is a special case of Corollary~\ref{cor:fn-vald} in which $p = 2$, the first continuous dissection $\F_1$ is trivial and $\F_2 = \F$ (see the second example following Definition~\ref{dfn:layered}).

It remains to prove continuity.  Let us write elements of $\C_\cts(Y,EA)$ as functions on $Y\times I$. Fix $x\in X$, and consider a neighbourhood of the identity in $\C_\cts(Y,EA)$ of the form
\[W := \{g:\ g(y,\cdot) \in V\ \forall y \in K\},\]
where $V$ is a neighbourhood of the identity in $EA$ and $K\subseteq Y$ is compact.  We must find a neighbourhood $U$ of $x$ in $X$ such that
\[f(x_1,\cdot,\cdot) - f(x,\cdot,\cdot) \in W \quad \forall x_1 \in U.\]
This will complete the proof, because such sets $W$ for different choices of $K$ and $V$ form a neighbourhood basis at the identity in the compact-open topology of $\C_\cts(Y,EA)$.

Since $K$ is compact and $\F$ is locally finite, $x$ has a neighbourhood $U_1$ such that $\F|_{U_1\times K}$ is finite, say of cardinality $m$. It follows that $F(x_1)|_{K\times I}$ lies in $\C_\cts(K,E^{(m)}A) \subseteq \C_\cts(K,EA)$ for all $x_1 \in U_1$, recalling that $E^{(m)}A$ is the set of members of $EA$ that have at most $m$ discontinuities.  Having found this $m$, there are an $\eps > 0$ and an identity neighbourhood $B \subseteq A$ such that
\[\big\{f \in E^{(m)}A:\ \l\{t:\ f(t) \in B\} > 1 - \eps\big\} \subseteq V\]
(observe that $\{t:\ f(t) \in B\}$ is a finite union of intervals, so certainly measurable).

However, again using the compactness of $K$, we may now find a possibly smaller neighbourhood $U \subseteq U_1$ such that the following two conditions hold:
\begin{itemize}
\item $|\xi(x_1,y) - \xi(x,y)| < (\eps/2m)$ for all $x_1 \in U$, $y \in K$ and $\xi \in \F$;
\item if $C$ is a minimal $\F$-wedge such that
\[C\cap (\{(x,y)\}\times I) \neq \emptyset \quad \hbox{and} \quad C\cap (\{(x_1,y)\}\times I) \neq \emptyset\]
for some $y \in K$ and $x_1 \in U$, and if $f_C:\ol{\b(C)}\to A$ is the corresponding continuous function promised by Definition~\ref{dfn:layered}, then
\[f_C(x,y) - f_C(x_1,y) \in B.\]
\end{itemize}

For each $y \in Y$, the interval $\{y\}\times I$ is partitioned into minimal subintervals of the form $(\xi(x,y),\xi'(x,y)]$ for certain pairs $\xi,\xi' \in \F$.  Each of these minimal subintervals describes the intersection of $\{(x,y)\}\times I$ with some minimal $\F$-wedge $C$. By the first condition above we also know that the end-points of the corresponding interval $(\xi(x_1,y),\xi'(x_1,y)]$ above $(x_1,y)$ are different from those of $(\xi(x,y),\xi'(x,y)]$ by less than $(\eps/2m)$ for any $x_1 \in U$.  Therefore, for any
\[t \in T := I\Big\backslash \bigcup_{\xi \in \F}\big(\xi(x,y) - \eps/2m,\ \xi(x,y) + \eps/2m\big)\]
and any $x_1 \in U$, the triples $(x,y,t)$ and $(x_1,y,t)$ lie in the same minimal $\F$-wedge $C$, and hence
\[f(x,y,t) - f(x_1,y,t) = f_C(x,y) - f_C(x_1,y)\in B,\]
using the second condition above.  Since the complement of $T$ is a union of at most $m$ intervals of length less than $\eps/m$, we also have $\l(T) > 1 - \eps$, and so the proof is complete. \qed

\begin{lem}\label{lem:closures-of-projections}
Let $C^- \subseteq X_{\leq p-1}\times I^{p-1}$ be a minimal $(\F_1,\ldots,\F_{p-1})$-multiwedge, $C_p \subseteq X_{\leq p}\times I$ be a minimal $\F_p$-wedge, and let $C$ be the resulting $(\F_1,\ldots,\F_p)$-multiwedge:
\[C = C^-\times_{X_{\leq p}} C_p.\]
Also, let
\begin{multline*}
\b^-:X_{\leq p-1}\times I^{p-1}\to X_{\leq p-1}, \quad \b:X_{\leq p}\times I^p\to X_{\leq p}\\ \hbox{and} \quad \k:X_{\leq p}\times I\to X_{\leq p}
\end{multline*}
be the coordinate projections. If $x^- \in \ol{\b^-(C^-)}$ and $(x^-,x_p) \in \k(C_p)$, then $(x^-,x_p) \in \ol{\b(C)}$.
\end{lem}

\noindent\textbf{Proof}\quad If
\[C_p = \{(y^-,y_p,t):\ \xi_1(y^-,y_p) < t \leq \xi_2(y^-,y_p)\},\]
then
\[\k(C_p) = \{(y^-,y_p):\ \xi_2(y^-,y_p) > \xi_1(y^-,y_p)\},\]
so this is an open set.  Therefore for any sufficiently small neighbourhood $U$ of $x^-$ one has $U\times \{x_p\} \subseteq \k(C_p)$, meaning that for any $y^- \in U$ there is some $t \in I$ such that $(y^-,x_p,t) \in C_p$. On the other hand, $U\cap \b^-(C^-) \neq \emptyset$ for any open set $U$ containing $x^-$, meaning that for some $y^- \in U$ and $t^- \in I^{n-1}$ one has $(y^-,t^-) \in C^-$.  Putting these together gives
\[(y^-,x_p,t^-,t) \in C \quad \hbox{and hence} \quad (U\times \{x_p\})\cap \b(C) \neq \emptyset.\]
Since $U$ was an arbitrarily small neighbourhood of $x^-$, this implies $(x^-,x_p) \in \ol{\b(C)}$, as required. \qed

\vspace{7pt}

\noindent\textbf{Proof of Proposition~\ref{prop:fn-vald-semi-layered}}\quad Corollary~\ref{cor:fn-vald} tells us that $F$ takes values in $\C_\cts(X_p,EA)$, so it remains to prove that it is semi-layered. Let $\b^-$, $\b$ and $\k$ be as in the previous lemma. We must show that for any minimal $(\F_1,\ldots,\F_{p-1})$-multiwedge
\[C^- := C_1 \times_{X_{\leq p-1}} \cdots \times_{X_{\leq p-1}} C_{p-1},\]
there is a continuous function $F_{C^-}:\ol{\b^-(C^-)}\to \C_\cts(X_p,EA)$ satisfying $F|_{C^-} = F_{C^-}\circ \b^-|_{C^-}$.

As in the proof of Corollary~\ref{cor:fn-vald}, if $C_p$ is a minimal $\F_p$-wedge and we write
\[C := C_1 \times_{X_{\leq p}} \cdots \times_{X_{\leq p}} C_{p-1}\times_{X_{\leq p}}C_p,\]
then there is a continuous function $f_{C_p}:\ol{\b(C)}\to A$ (indexing here by $C_p$ instead of $C$, since $C_1$, \ldots, $C_{p-1}$ are fixed) such that
\[F(x^-,t^-)(x_p,t) = f_{C_p}(x^-,x_p)\]
whenever $(x^-,t^-) \in C^-$ and $(x^-,x_p,t)\in C_p$. This already shows that for each $C^-$ the restriction $F|_{C^-}$ depends only on $x^-$, not on $t^-$.  It therefore defines a function $F_{C^-}:\b^-(C^-)\to \C_\cts(X_p,EA)$.  Moreover, by Lemma~\ref{lem:closures-of-projections} we may actually define $F_{C^-}(x^-)$ for any $x^- \in \ol{\b^-(C^-)}$ by
\[F_{C^-}(x^-)(x_p,t) = f_{C_p}(x^-,x_p) \quad \hbox{whenever}\ (x^-,x_p,t) \in C_p,\]
since $(x^-,x_p) \in \ol{\b(C)} = \rm{dom}(f_{C_p})$ whenever $(x^-,x_p) \in \k(C_p)$.  In these terms, we have just shown that
\[F|_{C^-} = F_{C^-}\circ\b^-|_{C^-}.\]

The proof is completed by showing that this $F_{C^-}$ is continuous.  To see this, define
\[f'_{C^-}:\ol{\b^-(C^-)}\times X_p\times I\to A\]
by the requirement that
\[f'_{C^-}(x^-,x_p,t) = f_{C_p}(x^-,x_p) \quad \hbox{whenever}\ (x^-,x_p,t) \in C_p.\]
This is manifestly a semi-layered function, semi-controlled by $\F_p$, and now $F_{C^-}$ is the function defined from $f'_{C^-}$ as in the statement of Lemma~\ref{lem:pre-fn-vald-semi-layered}.  That lemma therefore completes the proof.\qed

\section{Almost layered functions}\label{sec:al}

Now assume further that $A$ is a Polish topological group with a translation-invariant complete metric $d$.  In this setting another class of functions will come into play.

\begin{dfn}[Almost layered functions]\label{dfn:al}
A function $X_{\leq p}\times I^p\to A$ is \textbf{almost layered} if it is a uniform limit of layered functions.
\end{dfn}

\noindent\textbf{Remark}\quad It is important that this time we insist on a uniform limit, rather than just locally uniform. \fin

\vspace{7pt}

Like Definition~\ref{dfn:layered}, this implicitly makes reference to the structure of $X_{\leq p}$ as a product of $p$ spaces.

\begin{lem}\label{lem:sums-still-a-l}
If a function is a uniform limit of almost layered functions, then it is almost layered, and the sum of two almost layered functions is almost layered.
\end{lem}

\noindent\textbf{Proof}\quad The first part follows by the usual diagonal argument, and the second by a simple appeal to Lemma~\ref{lem:sums-still-s-l}. \qed

\vspace{7pt}

The following analog of Lemma~\ref{lem:pre-pullback} is also immediate, simply by pulling back layered approximants and applying Lemma~\ref{lem:pre-pullback} itself to those.

\begin{lem}[Pulling back and slicing]\label{lem:pullback}
Suppose that $\phi_i:X_{\leq i}\to Y_i$ is an ascending tuple of maps between metrizable spaces and that $f:Y_{\leq p}\times I^p\to A$ is an almost layered function.  Abbreviate $\phi_{\leq p} =: \phi$. Then the pullback $\phi^\ast f := f(\phi(\cdot),\cdot)$ is an almost layered function on $X_{\leq p}\times I^p$. \qed
\end{lem}

Analogously to semi-layered functions, almost layered functions can be lifted through quotients of Polish modules.  This proof is rather different from Proposition~\ref{prop:lift-s-l}, but is very similar to the proof of Proposition~\ref{prop:lifting}.

\begin{prop}[Lifting almost layered functions]\label{prop:lift-a-l}
Suppose that $B \into A \onto A/B$ is an exact sequence of Polish groups (but with no assumption of a continuous cross-section).  Then any almost layered function $f:X_{\leq p}\times I^p \to A/B$ has an almost layered lift $X_{\leq p}\times I^p\to A$.
\end{prop}

\noindent\textbf{Proof}\quad Consider $A/B$ endowed with the quotient $\ol{d}$ of the metric $d$. Let $d_\infty$ and $\bar{d}_\infty$ denote respectively the uniform metrics on spaces of $A$- and $(A/B)$-valued functions.

Let $(\g_m)_{m\geq 1}$ be a sequence of layered functions $X_{\leq p}\times I^p \to A/B$ such that $\bar{d}_\infty(f,\g_m) < 2^{-m}$, and for each $m$ let $(\F_{m,1},\ldots,\F_{m,p})$ be a tuple of l-complete continuous dissections that controls $\g_m$. We may assume that $\F_{m+1,i}\supseteq \F_{m,i}$ for each $m$ and $i$, for otherwise this can be arranged by replacing each $\F_{m,i}$ with $\F'_{m,i} := \ol{\bigcup_{m'\leq m}\F_{m',i}}$.

For each $m$ let $\P^0_m$ be the partition of $X_{\leq p}\times I^p$ into the level sets of $\g_m$, and let $\P_m := \bigvee_{m'\leq m}\P_{m'}^0$ (the common refinement).  Because the $\F_{m,i}$ are l-complete and non-decreasing in $m$, any cell $C \in \P_m$ is a union of $(\F_{m,1},\ldots,\F_{m,p})$-multiwedges.

Choose a layered lift $\hat{\g}_m$ of each $\g_m$ recursively as follows.  When $m=1$, for each $C \in \P_1$ we simply choose a lift $\hat{\g}_m(C) \in A$ of $\g_m(C) \in A/B$.  Now suppose we have already constructed $\hat{\g}_m$ for some $m$.  Then each $C \in \P_{m+1}$ is contained in some $C_0 \in \P_m$, and picking a reference point $(x,t) \in C$ we know that
\[\bar{d}(\g_{m+1}(C),\g_m(C_0)) \leq \bar{d}(\g_{m+1}(C),f(x,t)) + \bar{d}(f(x,t),\g_m(C_0)) < 2^{-m+1}.\]
By the definition of $\bar{d}$ as a quotient metric, this implies that there is some lift of the point $\g_{m+1}(C)$ lying within $d$-distance $2^{-m+2}$ of $\hat{\g}_m(C_0)$.  Define $\hat{\g}_{m+1}(C)$ to be such a lift.

Each $\hat{\g}_m$ is a lift of $\g_m$ which is layered and controlled by $(\F_{m,1},\ldots,\F_{m,p})$, and the sequence of functions $(\hat{\g}_m)_{m\geq 1}$ is uniformly Cauchy. Letting $F$ be its uniform limit gives an almost layered lift of $f$. \qed

\vspace{7pt}

The next lemma shows that the definition of almost layered functions is insensitive to enlargement of the target module. 

\begin{lem}\label{lem:layered-in-subgroup}
If $B$ is a Polish group, $A$ is a closed subgroup and $f:X_{\leq p}\times I^p\to A$ is almost layered as a $B$-valued function, then it is almost layered as an $A$-valued function.
\end{lem}

\noindent\textbf{Proof}\quad Suppose that $\eps > 0$ and let $\g:X_{\leq p}\times I^p\to B$ be a layered function satisfying $d_\infty(f,\g) < \eps$.  Let $\P$ be the level-set partition of $\g$.  Then for every $C \in \P$, the single value $\g(C)$ must lie within $\eps$ of all the values taken by $f$ on $C$. Defining $\g':X_{\leq p}\times I^p\to A$ to take a constant value lying in $f(C)$ for each such $C$ therefore gives a new layered function which is $A$-valued and satisfies $d_\infty(f,\g') < 2\eps$.  Since $\eps$ was arbitrary this completes the proof. \qed

\vspace{7pt}

It is clear that any almost layered function is measurable.  The following result provides the link between semi-layered and almost-layered functions.

\begin{lem}\label{lem:cts-implies-a-l}
If $f:X_{\leq p}\times I^p \to A$ is a semi-layered function, say semi-controlled by $(\F_1,\ldots,\F_p)$, then $f$ is almost layered.
\end{lem}

\noindent\textbf{Proof}\quad Given $\eps > 0$ we must find a layered function that is uniformly $\eps$-close to $f$.

Let $\P$ be the partition of $X_{\leq p}\times I^p$ into minimal $(\F_1,\ldots,\F_p)$-multiwedges, and as previously let $\b:X_{\leq p}\times I^p\to X_{\leq p}$ be the coordinate projection. For each $C \in \P$, let $f_C$ be a continuous function on $\ol{\b(C)}$ such that $f|_C = f_C\circ \b|_C$.  By continuity, each $x \in \ol{\b(C)}$ has a neighbourhood $W_{C,x}$ such that $f_C(\ol{\b(C)}\cap W_{C,x})$ lies within the $(\eps/2)$-ball around $f_C(x)$.  Moreover, since $x$ can lie in $\ol{\b(C)}$ for only finitely many $C \in \P$, the resulting intersection $U_x := \bigcap_{C:\,\ol{\b(C)}\ni x} W_{C,x}$ is still a neighbourhood of $x$.

The collection $\cal{U}$ of these $U_x$ is an open cover of $X_{\leq p}$.  Therefore Lemma~\ref{lem:wedges-in-cover} promises an l-complete continuous dissection $\cal{G}$ whose minimal wedges are all contained in $\b$-pre-images of elements of $\cal{U}$.

Let $\F := \ol{\F_p\cup \cal{G}}$, and consider a minimal $(\F_1,\ldots,\F)$-multiwedge $D$.  Since $\F \supseteq \F_p$, $D$ is wholly contained in some minimal $(\F_1,\ldots,\F_p)$-multiwedge, say $C$. Since also $\F\supseteq \cal{G}$, $D$ is also contained in some set of the form $U_x\times I^p \subseteq W_{C,x}\times I^p$. By the construction of the sets $W_{C,x}$, this implies that the image $f(D)$ has $d$-diameter at most $\eps$.  Thus we obtain a layered function $\g$ which is $\eps$-close to $f$ by letting $\g$ take a fixed value from the image $f(D)$ for each such $D$. This completes the proof. \qed

\vspace{7pt}

Proposition~\ref{prop:fn-vald-semi-layered} quickly implies the following simple analog for almost layered functions.

\begin{lem}\label{lem:efface-al}
Suppose that each $X_i$ is a locally compact second countable metrizable space and that $X_p$ carries a Radon probability measure $\mu$, and let $\C(X_p,LA)$ denote the Polish group of $\mu$-equivalence classes of measurable maps $X_p\to LA$ with the topology of convergence in measure on compact sets.

If $f:X_{\leq p}\times I^p\to A$ is almost layered then the map $F:X_{\leq p-1}\times I^{p-1}\to A^{X_p\times I}$ defined by
\[F(x_1,\ldots,x_{p-1},t_1,\ldots,t_{p-1}) := f(x_1,\ldots,x_{p-1},\cdot,t_1,\ldots,t_{p-1},\cdot)\]
takes values in $\C(X_p,LA)$ and is almost layered for that target module.
\end{lem}

\noindent\textbf{Proof}\quad Let $\g_m$ be a sequence of layered functions such that $d_\infty(f,\g_m) < 2^{-m}$, and for each $m$ let
\[\eta_m(x_1,\ldots,x_{p-1},t_1,\ldots,t_{p-1}) := \g_m(x_1,\ldots,x_{p-1},\cdot,t_1,\ldots,t_{p-1},\cdot).\]
Then each $\eta_m$ defines a semi-layered function to $\C_\cts(X_p,EA)$ by Proposition~\ref{prop:fn-vald-semi-layered}, and hence also to $\C(X_p,LA)$ (since the obvious homomorphisms
\[\C_\cts(X_p,EA) \to \C_\cts(X_p,LA)\to \C(X_p,LA)\]
are both continuous).  Moreover, for each $x_1$, \ldots, $x_{p-1}$, $t_1$, \ldots, $t_{p-1}$ we have
\[d_\infty(\eta_m(x_1,\ldots,x_{p-1},t_1,\ldots,t_{p-1}), F(x_1,\ldots,x_{p-1},t_1,\ldots,t_{p-1})) < 2^{-m}\]
as $m\to\infty$, where $d_\infty$ denotes the uniform metric on functions $X_p\times I\to A$.  This is certainly stronger than the topology on $\C(X_p,LA)$, so this shows that $\eta_m$ converges uniformly to $F$ among functions $X_{p-1}\times I^{p-1}\to \C(X_p,LA)$.  Hence the proof is complete by Lemmas~\ref{lem:cts-implies-a-l} and~\ref{lem:sums-still-a-l}. \qed

\vspace{7pt}

Before turning to applications, we prove one more technical property of almost layered functions that will be crucial later.

\begin{lem}\label{lem:loc-to-glob-approx}
Suppose that $f:X_{\leq p}\times I^p\to A$ is a function with the property that for every $\eps > 0$ and every $x_1 \in X_1$ there are a neighbourhood $U$ of $x_1$ and a semi-layered function $\g_U:X_{\leq p}\times I^p\to A$ such that
\begin{multline}\label{eq:local-approx}
d\big(f(x_1,\ldots,x_p,t_1,\ldots,t_p),\g_U(x_1,\ldots,x_p,t_1,\ldots,t_p)\big) < \eps \\ \forall (x_1,\ldots,x_p,t_1,\ldots,t_p) \in U\times X_2\times \cdots \times X_p\times I^p.
\end{multline}
Then $f$ is almost layered.
\end{lem}

\noindent\textbf{Remark}\quad Heuristically, this lemma allows us to `localize' the approximability by layered functions without changing the class of almost layered functions, provided that localization is only in the first coordinate of $X_{\leq p}$. \fin

\vspace{7pt}

\noindent\textbf{Proof}\quad Let $\cal{U}$ be the open cover of $X_1$ by the sets appearing in the hypotheses, and for each $U \in \cal{U}$ let $(\F_{U,1},\ldots,\F_{U,p})$ be an ascending tuple of l-complete continuous dissections that controls $\g_U$.  

From these data, Lemma~\ref{lem:ctsdiss-loc-to-glob-approx} gives another l-complete ascending tuple $(\F_1,\ldots,\F_p)$ such that for every minimal $(\F_1,\ldots,\F_p)$-multiwedge $C$ there is some $U_C \in \cal{U}$ such that
\begin{itemize}
\item $C \subseteq U_C\times X_2\times \cdots \times X_p \times I^p$, and
\item $C$ is contained in some $(\F_{U_C,1},\ldots,\F_{U_C,p})$-multiwedge.
\end{itemize}

Now define $\g:X_{\leq p}\times I^p\to A$ by the stipulation that on each such $C$ it agrees with $\g_{U_C}$.  This is well-defined by the second property above, and it manifestly gives another semi-layered function.  Moreover, by the first property above and the assumed approximation of $f$ by $\g_U$ on $U\times X_2\times \cdots\times X_p\times I^p$, we now have $d(f,\g) < \eps$ everywhere.  Lemmas~\ref{lem:cts-implies-a-l} and~\ref{lem:sums-still-a-l} complete the proof. \qed

\section{Comparison of cohomology theories}\label{sec:A}

We can now prove the two key results that will give us comparable cocycle representations for $\rmH^\ast_\m$ and $\rmH^\ast_\Seg$.

The first fact we need is the following.

\begin{lem}\label{lem:d-a-l-is-a-l}
If $A$ is any Hausdorff topological group and $\s:G^{p+1}\times I^{p+1}\to A$ is a semi-layered cochain, then $d\s:G^{p+2}\times I^{p+2}\to A$ is also semi-layered.  If $A$ is Polish then the analogous fact holds among almost layered functions.
\end{lem}

\noindent\textbf{Proof}\quad In view of the defining formula
\begin{multline*}
d\s(g_1,\ldots,g_{p+2},t_1,\ldots,t_{p+2})\\ = \sum_{i=1}^{p+2}(-1)^{p+2 - i}\s(g_1,\ldots,\hat{g_i},\ldots,g_{p+2},t_1,\ldots,\hat{t_i},\ldots,t_{p+2}),
\end{multline*}
this follows at once from Lemmas~\ref{lem:pre-pullback} and~\ref{lem:sums-still-s-l} (in the semi-layered case) and Lemmas~\ref{lem:pullback},~\ref{lem:sums-still-a-l} and~\ref{lem:layered-in-subgroup} (in the almost layered case). \qed

\vspace{7pt}

Now, if $G$ is a metrizable topological group and $A$ is a $G$-module in Segal's category, let $\C^p_\rm{sl}(G,A)$ be the Abelian group of all $G$-equivariant semi-layered functions $G^{p+1}\times I^{p+1}\to A$.  Using these, form the complex
\[0 \to \C^1_\sl(G,A)\to \C_\sl^2(G,A) \to \cdots\]
with the alternating-sum differentials, which is well-defined by Lemma~\ref{lem:d-a-l-is-a-l}.  Finally, define $\rmH^\ast_\sl(G,A)$ to be the homology of this complex, and call this the \textbf{semi-layered cohomology of $(G,A)$}.

Similarly, if $G$ is l.c.s.c. and $A$ is a Polish $G$-module, let $\C_\al^p(G,A)$ denote the $G$-equivariant almost layered functions $G^{p+1}\times I^{p+1}\to A$, and form the complex
\[0 \to \C^1_\al(G,A)\to \C_\al^2(G,A) \to \cdots\]
with the alternating-sum differentials.  Let $\rmH^\ast_\al(G,A)$ be its homology, and call this the \textbf{almost layered cohomology of $(G,A)$}.  It is worth emphasizing that while elements of $\C^p_\al(G,A)$ are equivariant, it may not be possible to find layered functions that approximate them and are equivariant.

\begin{prop}\label{prop:Seg=sl}
If $G$ is a topological group in the category of k-spaces, then $\rmH^\ast_\sl$ defines a connected sequence of functors on Segal's category of $G$-modules which is isomorphic to $\rmH^\ast_\Seg$.
\end{prop}

\begin{prop}\label{prop:m=al}
If $G$ is l.c.s.c., then $\rmH^\ast_\al$ defines a connected sequence of functors on Polish $G$-modules which is isomorphic to $\rmH^\ast_\m$.
\end{prop}

Both of these propositions will be proved via Buchsbaum's criterion.  In each case we must check (i) the degree-zero interpretation, (ii) the construction of a long exact sequence and (iii) effaceability on the relevant category of modules.  The switchback maps of the long exact sequence will be constructed in the process.  These will be fairly simple consequences of the properties of semi- and almost layered functions established in the previous sections.  However, let us first complete proof of our main result.

\vspace{7pt}

\noindent\textbf{Proof of Theorem A from Propositions~\ref{prop:Seg=sl} and~\ref{prop:m=al}}\quad If $A$ is discrete, then any uniformly convergent sequence of $A$-valued functions must stabilize after finitely many terms, so in this setting semi-layered and almost layered functions are all actually just layered.  Hence the defining complexes of $\rmH^\ast_\rm{sl}(G,A)$ and $\rmH^\ast_\rm{al}(G,A)$ are the same. \qed

\vspace{7pt}

For a general Polish module $A$ which is locally contractible, Lemma~\ref{lem:cts-implies-a-l} gives a comparison map
\[\rmH^\ast_\Seg \cong \rmH^\ast_\sl \to \rmH^\ast_\al \cong \rmH^\ast_\m,\]
but it seems unlikely that it is always an isomorphism (see also the results of~\cite{AusMoo--cohomcty}).

\subsection*{Segal and semi-layered theories}

Most of the remaining work for the semi-layered theory is in establishing the long exact sequence.  This will need an analog of Proposition~\ref{prop:lift-s-l} for equivariant functions.

\begin{lem}\label{lem:equiv-sl-lift}
Suppose that $B\into A\onto A/B$ is an exact sequence of Hausdorff topological Abelian groups that admits a local continuous cross-section.  Then any equivariant semi-layered function $f:G^{p+1}\times I^{p+1}\to A/B$ has an equivariant semi-layered lift $G^{p+1}\times I^{p+1} \to A$.
\end{lem}

\noindent\textbf{Proof}\quad This follows by combining Proposition~\ref{prop:lift-s-l} and Lemma~\ref{lem:pre-pullback}.  Suppose that $f$ is semi-controlled by $(\F_1,\ldots,\F_{p+1})$ with each $\F_i$ being l-complete.

Let $X_1 := \{e\}$ and $X_i := G$ for $2 \leq i \leq p+1$; clearly these are still metrizable topological spaces.  Applying Lemma~\ref{lem:pre-pullback} to the ascending tuple of maps $\phi_i:X_{\leq i}\to G^i$ defined by
\[\phi_1(e) = e \quad \hbox{and} \quad \phi_i(e,g_2,\ldots,g_i) = g_i\ \hbox{for}\ i\geq 2,\]
we find that the restriction $f|_{\{e\}\times G^p\times I^{p+1}}$ is semi-layered and semi-controlled by $(\phi_1^\ast\F_1,\ldots,\phi_{\leq p+1}^\ast\F_{p+1})$.

Therefore, applying Proposition~\ref{prop:lift-s-l} to this restriction gives a semi-layered lift $F_0:\{e\}\times G^p\times I^{p+1}\to A$.  Suppose that $F_0$ is semi-controlled by the tuple $(\cal{G}_1,\ldots,\cal{G}_{p+1})$.

Lastly, let $F:G^{p+1}\times I^{p+1}\to A$ be the extension of $F_0$ determined by equivariance:
\[F(g_1,\ldots,g_{p+1},t_1,\ldots,t_{p+1}) = g_1\cdot\big(F_0(e,g_1^{-1}g_2,\ldots,g_1^{-1}g_{p+1},t_1,\ldots,t_{p+1})\big).\]
Since $f$ was equivariant, $F$ must be a lift of $f$.  We will show in two further steps that $F$ is also semi-layered.

Define $F_1$ by
\[F_1(g_1,\ldots,g_{p+1},t_1,\ldots,g_{p+1}) = F_0(e,g_1^{-1}g_2,\ldots,g_1^{-1}g_{p+1},t_1,\ldots,t_{p+1}).\]
Then this is equal to $\psi_{\leq p+1}^\ast F_0$, where $\psi_{\leq p+1}:G^{p+1}\to \{e\}\times G^p$ is obtained from the ascending tuple of functions
\[\psi_i:G^i\to X_i:(g_1,\ldots,g_i) \to g_1^{-1}g_i.\]
Therefore Lemma~\ref{lem:pre-pullback} shows that $F_1$ is semi-layered, semi-controlled by $(\psi_1^\ast\cal{G}_1,\ldots,\psi_{\leq p+1}^\ast\cal{G}_{p+1})$.  Let $\cal{P}$ be the partition of $G^{p+1}\times I^{p+1}$ into minimal $(\psi_1^\ast\cal{G}_1,\ldots,\psi_{\leq p+1}^\ast\cal{G}_{p+1})$-multiwedges.

Observe that $F(g_1,\ldots) := g_1\cdot(F_1(g_1,\ldots))$.  We will prove that $F$ also satisfies Definition~\ref{dfn:layered} with the same partition $\cal{P}$.  As previously, let $\b:G^{p+1}\times I^{p+1}\to G^{p+1}$ be the coordinate projection.  If $C \in \cal{P}$, then there is a continuous function $h_C:\ol{\b(C)}\to A$ such that $F_1|_C = h_C\circ \b|_C$.  For $(g_1,\ldots,g_{p+1},t_1,\ldots,t_{p+1})\in C$ this now gives
\[F(g_1,\ldots,g_{p+1},t_1,\ldots,t_{p+1}) = g_1\cdot(h_C(g_1,\ldots,g_{p+1})),\]
so defining $h'_C(g_1,\ldots,g_{p+1}) := g_1(h_C(g_1,\ldots,g_{p+1}))$, this is also a continuous function on $\ol{\b(C)}$ whose lift gives the restriction $F|_C$.  This completes the proof. \qed

\vspace{7pt}

\begin{cor}
The theory $\rmH^\ast_\sl(G,\cdot)$ has long exact sequences on Segal's category.
\end{cor}

\noindent\textbf{Proof}\quad This follows the standard pattern.  Suppose that $B\into A\onto A/B$ is an exact sequence of modules.  Then the switchback maps $\rmH^p_\sl(G,A/B)\to \rmH^{p+1}_\sl(G,B)$ are defined cocycle-wise.  If $\s:G^{p+1}\times I^{p+1}\to A/B$ is a semi-layered cocycle, Lemma~\ref{lem:equiv-sl-lift} gives an equivariant semi-layered lift $\tau:G^{p+1}\times I^{p+1}\to A$, whose coboundary $d\tau$ must take values in $B$ because $d\s = 0$.  The image of $[\s]$ under the switchback is defined to be $[d\tau]$; this is well-defined because if $\s$ were a semi-layered coboundary, say $\s = d\a$, then another appeal to Lemma~\ref{lem:equiv-sl-lift} gives an equivariant semi-layered lift of $\a$, say $\b$, and hence
\[\tau = d\b + (B\hbox{-valued}) \quad \Rightarrow \quad d\tau = d(B\hbox{-valued}),\]
so $[d\tau] = 0$.

The remaining step is to verify that the resulting sequence
\begin{multline*}
\ldots \to \rmH^p(G,B) \to \rmH^p(G,A)\to \rmH^p(G,A/B)\\
\stackrel{\rm{switchback}}{\to}\rmH^{p+1}(G,B) \to \ldots
\end{multline*}
is exact; this follows just as in classical discrete group cohomology, since Lemma~\ref{lem:equiv-sl-lift} guarantees that lifts may be chosen to be semi-layered wherever necessary. \qed

\vspace{7pt}

\noindent\textbf{Proof of Theorem~\ref{prop:Seg=sl}}\quad We check the three axioms in turn.

In degree zero, there are no semi-layered coboundaries, and a semi-layered cocycle is a semi-layered map $f:G\times I\to A$ such that, on the one hand, $f(g,t) - f(g',t') = 0$, so $f$ is constant, and on the other $f$ is equivariant, so that its constant value must lie in $A^G$.

Next we prove effaceability.  If $\s:G^{p+1}\times I^{p+1}\to A$ is a semi-layered cocycle semi-controlled by $(\F_1,\ldots,\F_{p+1})$, then setting
\[F(g_1,\ldots,g_p,t_1,\ldots,t_p)(g,t) := \s(g_1,\ldots,g_p,g,t_1,\ldots,t_p,t)\]
defines a map $G^p\times I^p\to A^{G\times I}$.  By Proposition~\ref{prop:fn-vald-semi-layered}, it takes values in $\C_\cts(G,EA)$, and when that module is given Segal's topology this map is semi-layered and semi-controlled by $(\F_1,\ldots,\F_p)$. Lastly, the $G$-equivariance of $F$ follows immediately from that of $\s$. Therefore Segal's embedding $A \into \C_\cts(G,EA)$ effaces semi-layered cohomology, just as it does $\rmH^\ast_\Seg$: the coboundary of the new cochain $F$ is equal to $\s$ by the same calculation as in the discrete-groups case.

Lastly, the long exact sequence has been constructed in the previous corollary. \qed

\subsection*{Measurable and almost layered theories}

Now we need analogous results for almost layered cohomology.

\begin{lem}[Lifting almost layered cocycles]\label{lem:equiv-lift-a-l}
If $B \into A \onto A/B$ is an exact sequence of Polish Abelian groups, then any $G$-equivariant almost layered function $f:G^{p+1}\times I^{p+1} \to A/B$ has an almost layered lift $G^{p+1}\times I^{p+1}\to A$.
\end{lem}

\noindent\textbf{Proof}\quad This mostly follows the same pattern as Lemma~\ref{lem:equiv-sl-lift}: this time we combine Proposition~\ref{prop:lift-a-l} and Lemma~\ref{lem:pullback}.

If $f:G^{p+1}\times I^{p+1}\to A$ is equivariant and almost layered, then applying Lemma~\ref{lem:pullback} to the maps $\phi_i:X_{\leq i}\to G^i$ defined by
\[\phi_1(e) = e \quad \hbox{and} \quad \phi_i(e,g_2,\ldots,g_i) = g_i\ \hbox{for}\ i\geq 2\]
gives that the restriction $f|_{\{e\}\times G^p\times I^{p+1}}$ is almost layered.  Proposition~\ref{prop:lift-a-l} therefore gives gives a semi-layered lift $F_0:\{e\}\times G^p\times I^{p+1}\to A$ of this restriction. Now let $F:G^{p+1}\times I^{p+1}\to A$ be the extension of $F_0$ determined by equivariance:
\[F(g_1,\ldots,g_{p+1},t_1,\ldots,t_{p+1}) = g_1\cdot\big(F_0(e,g_1^{-1}g_2,\ldots,g_1^{-1}g_{p+1},t_1,\ldots,t_{p+1})\big).\]
Since $f$ was equivariant, $F$ must be a lift of $f$.  We will show in two further steps that $F$ is also almost layered.

First, the function
\[F_1(g_1,\ldots,g_{p+1},t_1,\ldots,t_{p+1}) := F_0(e,g_1^{-1}g_2,\ldots,g_1^{-1}g_{p+1},t_1,\ldots,t_{p+1})\]
is equal to $\psi_{\leq p+1}^\ast F_0$, where $\psi_{\leq p+1}:G^{p+1}\to \{e\}\times G^p$ is as in Lemma~\ref{lem:equiv-sl-lift}. Lemma~\ref{lem:pullback} shows that $F_1$ is almost layered.  We may therefore choose a sequence of layered functions $\g_m:G^{p+1}\times I^{p+1}\to A$ that converge to $F_1$ uniformly.

Consider the functions
\[\g'_m(g_1,\ldots,g_{p+1},t_1,\ldots,t_{p+1}) := g_1\cdot\big(\g_m(g_1,\ldots,g_{p+1},t_1,\ldots,t_{p+1})\big).\]
By the continuity of the $G$-action, for each $g_1 \in G$ there is an $\eps_{g_1} > 0$ such that
\[d(g_1x,g_1y) = d(g_1(x-y),0) < \eps/2 \quad \hbox{whenever}\ d(x,y) < \eps_{g_1},\]
and knowing this, another appeal to continuity gives a neighbourhood of the identity $W$ in $G$ such that 
\[d(gx,gy) = d\big((gg_1^{-1})g_1(x-y),0\big) < \eps \quad \hbox{whenever}\ g \in Wg_1\ \hbox{and}\ d(x,y) < \eps_{g_1}.\]

The sets $Wg_1$, $g_1 \in G$, form a cover, so since $G$ is metrizable we may choose a locally finite subcover $\cal{U}$.  Since each $U \in \cal{U}$ is contained in some $Wg_1$, the above inequality gives some $\eps_U > 0$ such that
\[d(gx,gy) < \eps \quad \hbox{whenever}\ g \in U\ \hbox{and}\ d(x,y) < \eps_U.\]

Now choose for each $U$ some $m_U \geq 1$ such that $d_\infty(F_1,\g_{m_U}) < \eps_U$, and hence
\begin{eqnarray*}
&&d\big(F(g_1,\ldots,g_{p+1},t_1,\ldots,t_{p+1}),\g'_{m_U}(g_1,\ldots,g_{p+1},t_1,\ldots,t_{p+1})\big)\\ &&= d\big(g_1\cdot\big(F_1(g_1,\ldots,g_{p+1},t_1,\ldots,t_{p+1})\big),g_1\cdot\big(\g_{m_U}(g_1,\ldots,g_{p+1},t_1,\ldots,t_{p+1})\big)\big)\\ &&< \eps
\end{eqnarray*}
for any $(g_1,\ldots,g_{p+1},t_1,\ldots,t_{p+1}) \in U\times G^p\times I^{p+1}$.  This is the condition required by Lemma~\ref{lem:loc-to-glob-approx}, so $F$ is almost layered. \qed

\vspace{7pt}

\noindent\textbf{Proof of Proposition~\ref{prop:m=al}}\quad Once again this follows by Buchsbaum's criterion. For any $A$ the group $\rmH^0_\al(G,A)$ is identified with $A^G$ just as in the semi-layered case.  Effacement also follows as in the semi-layered case, this time using Lemma~\ref{lem:efface-al}.  Lastly, the long exact sequence follows by the standard construction using Lemma~\ref{lem:equiv-lift-a-l}.  \qed

\bibliographystyle{amsalpha}
\bibliography{bibfile}

\end{document}